\documentclass[dvips,twoside]{article}
\usepackage[a5paper,layoutwidth=148.5mm,layoutheight=7.77817in,%outer=28pt,inner=38pt,
hmargin=40bp,
vmargin=40bp
]{geometry}

\usepackage{microtype}
\usepackage[utf8]{inputenc}
\usepackage[T1]{fontenc}
\usepackage{lmodern}
\usepackage{fix-cm}
\usepackage{textcomp}
\usepackage{eufrak}
\usepackage[normalem]{ulem}
\renewcommand\uline[1]{\textbf{#1}}
\usepackage{soulutf8}
\sodef\so{}{.1em}{.5em plus.1em minus.1em}{.5em plus.1em minus.1em}

\usepackage{pbox}
\usepackage{ebproof}
\usepackage{xspace}
\usepackage{parallel}
\ebproofset{separation=1em,center=false}
\usepackage{mathtools}
\newcommand\Land{\mathbin{\&}}
\newcommand\rmM{\mathrm{M}}
\newcommand\rmP{\mathrm{P}}
\newcommand\rmK{\mathrm{K}}
\newcommand\rmA{\mathrm{A}}
\newcommand\rmB{\mathrm{B}}
\newcommand\rmN{\mathrm{N}}
\newcommand\itM{M}
\newcommand\red[1]{\so{#1}}
\newcommand\blue[1]{\emph{#1}}

\newcommand\bs{} \newcommand\nbsp{~}
\newcommand\bsp{\space}
\newcommand\bruch[1]{\ensuremath{\mid^{#1}}}

\usepackage[hyphens]{url}\RequirePackage{etoolbox}
\appto\UrlSpecials{\do\R{\penalty0 \mathchar`\R }\do\U{\penalty0 \mathchar`\U }\do\D{\penalty0 \mathchar`\D }\do\F{\penalty0 \mathchar`\F }\do\A{\penalty0 \mathchar`\A }\do\W{\penalty0 \mathchar`\W }\do\Z{\penalty0 \mathchar`\Z }\do\G{\penalty0 \mathchar`\G }\do\P{\penalty0 \mathchar`\P }\do\N{\penalty0 \mathchar`\N }}
\providecommand\url[1]{\texttt{\detokenize{#1}}}
\usepackage{natbib}
\let\bibfont=\small
\makeatletter\def\@biblabel#1{#1.}
\makeatother

\bibpunct{(}{)}{,}{a}{}{,}
\newcommand\urlprefix{}
\newcommand\bbletal{et al.}
\newcommand\bbland{and}
\newcommand\bbled{ed.}
\newcommand\bbleds{eds.}

\newcommand\bblin{in}

\newcommand\bblphdthesis{Ph.D. thesis}
\newcommand\natexlab[1]{\textit{#1}}
\newcommand{\enquote}[1]{`#1'}
\providecommand{\citenamefont}[1]{#1}
\newcommand{\Capitalize}[1]{\uppercase{#1}}
\newcommand{\capitalize}[1]{\expandafter\Capitalize#1}

\usepackage{amsfonts}

\usepackage[english,dutch,french,german,british]{babel}
\newcommand\german[1]{\foreignlanguage{german}{#1}}
\usepackage{lmodern}
\usepackage{fix-cm}

\usepackage[T1]{fontenc}
\usepackage[utf8]{inputenc}
\usepackage{amsmath}
\usepackage{amssymb}
\newcommand\bigland{\mathop{\textstyle\bigwedge}\nolimits}
\newcommand\lt{\leqslant}
\renewcommand*\neg[1]{\setbox0\hbox{$\mathaccent"0362{#1}^H$}\setbox2\hbox{$\mathaccent"0362{\kern0pt#1}^H$}\ifdim\ht0=\ht2 \overline{#1}\else \bar#1\fi}

\newcommand\Deutsch[1]{}

\usepackage{enumitem}
\PassOptionsToPackage{urlcolor=magenta,citecolor=blue,anchorcolor=blue,linkcolor=blue,menucolor=cyan}{hyperref}
\usepackage[pdftitle=,pdfauthor=,pdfdisplaydoctitle=true,pdflang=en-GB,bookmarksopen=false,breaklinks=true,unicode=true,hyperindex=true,pdfstartview=FitH,bookmarksnumbered=true]{hyperref}

\usepackage{doi}

\ifx\texorpdfstring\undefined\newcommand\texorpdfstring[2]{#1}\fi
\ifx\pdfbookmark[3]{}\fi

\usepackage{breakurl}
\usepackage{authblk}

\begin{document}
\title{\vspace*{-20bp}Lorenzen's proof of consistency for elementary number theory}
\author{Thierry Coquand\\Computer science and engineering department, University of Gothenburg, Sweden, \url{coquand@chalmers.se}.\and Stefan Neuwirth\\Laboratoire de mathématiques de Besançon, Université Bourgogne Franche-Comté, France, \url{stefan.neuwirth@univ-fcomte.fr}.%
}

\date{}
\maketitle

\begin{abstract}
  We present a manuscript of Paul Lorenzen that provides a proof of consistency for elementary number theory as an application of the construction of the free countably complete pseudocomplemented semilattice over a preordered set. This manuscript rests in the Oskar-Becker-Nachlass at the Philosophisches Archiv of Universität Konstanz, file OB 5-3b-5. It has probably been written between March and May 1944. We also compare this proof to Gentzen's and Novikov's, and provide a translation of the manuscript.

\end{abstract}

We present a manuscript of Paul Lorenzen that arguably dates back to 1944 and provide an edition and a translation, with the kind permission of Lorenzen's daughter, Jutta Reinhardt.

It provides a constructive proof of consistency for elementary number theory by showing that it is a part of a trivially consistent cut-free calculus. The proof resorts only to the inductive definition of formulas and theorems.

More precisely, Lorenzen proves the admissibility of cut by double induction, on the complexity of the cut formula and of the derivations, without using any ordinal assignment, contrary to the presentation of cut elimination in most standard texts on proof theory.

Prior to that, he proposes to define a countably complete pseudocomplemented semilattice as a deductive calculus, and shows how to present it for constructing the free countably complete pseudocomplemented semilattice over a given preordered set. 

He arrives at the understanding that the existence of this free kind of lattice captures the formal content of the consistency of elementary number theory, the more so as he has come to understand that the existence of another free kind of lattice captures the formal content of ideal theory. In this way, lattice theory provides a bridge between algebra and logic: by the concept of preorder, the divisibility of elements in a ring becomes commensurate with the material implication of numerical propositions; the lattice operations give rise to the ideal elements in algebra and to the compound propositions in logic.

The manuscript has remained unpublished, being superseded by Lorenzen's `Algebraische und logistische Untersuchungen über freie Verbände' that appeared in 1951 in \textit{The Journal of Symbolic Logic}. These `Algebraic and logistic investigations on free lattices' have immediately been recognised as a landmark in the history of infinitary proof theory, but their approach and method of proof have not been incorporated into the corpus of proof theory.

\section{The beginnings}\label{sec:preparation}

In 1938, Paul Lorenzen defends his Ph.D.\ thesis under the supervision of Helmut Hasse at Göttingen, an `Abstract foundation of the multiplicative ideal theory', i.e.\ a foundation of divisibility theory upon the theory of cancellative monoids. He is in a process of becoming more and more aware that lattice theory is the right framework for his research. \citet[footnote on p.~536]{lorenzen39} thinks of understanding a system of ideals as a lattice, with a reference to \citealt{koethe37}; in the definition of a semilattice-ordered monoid on p.~544, he credits Dedekind's two seminal articles of 1897 and 1900 for developing the concept of lattice. On 6 July 1938 he reports to Hasse: `Momentarily, I am at making a lattice-theoretic excerpt for Köthe'.\footnote{\german{Helmut-Hasse-Nachlass, Niedersächsische Staats- und Universitätsbibliothek Göttingen}, Cod.\ Ms.~H.~Hasse 1:1022, \citealp[edited in][§~4]{neuwirthkonstanz}.} He also reviews several articles on this subject for the \emph{Zentralblatt}, e.g.\ \citealt{klein39} and \citealt{george39} which both introduce semilattices, \citealt{whitman41} which studies free lattices. He also knows about the representation theorem for boolean algebras in \citealt{stone36} and he discusses the axioms for the arithmetic of real numbers in \citealt{tarski37} with Heinrich Scholz.\footnote{See the collection of documents grouped together by Scholz under the title `\german{\emph{Paul Lorenzen}: Gruppentheoretische Charakterisierung der reellen Zahlen} [Group theoretic characterisation of the real numbers]' and deposited at the \german{Bibliothek des Fachbereichs Mathematik und Informatik} of the \german{Westf{\"a}lische Wil\-helms-Uni\-ver\-si\-t{\"a}t M{\"u}nster}, as well as several letters filed in the \german{Heinrich-Scholz-Archiv} at \german{Universitäts- und Landesbibliothek Münster}, the earliest dated 7 April 1944.}

In 1939, he becomes assistant to \german{Wolfgang Krull} at Bonn. During World War~II, he serves first as a soldier and then, from 1942 on, as a teacher at the naval college \german{Wesermünde}. He devotes his `off-duty evenings all alone on [his] own'\footnote{Carbon copy of a letter to Krull, 13 March 1944, \german{Paul-Lorenzen-Nachlass, Philosophisches Archiv, Universität Konstanz}, PL~1-1-131,  \citealp[edited in][§~6]{neuwirthkonstanz}.} to mathematics with the goal of habilitating. On 25 April 1944, he writes to his advisor that `[\dots]\ it became clear to me---about 4 years ago---that a system of ideals is nothing but a semilattice'.\footnote{Carbon copy of a letter to Krull, PL~1-1-132, \citealp[edited in][§~6]{neuwirthkonstanz}.}

He will later recall a talk by Gerhard Gentzen on the consistency of elementary number theory in 1937 or 1938 as a trigger for his discovery that the reformulation of ideal theory in lattice-theoretic terms reveals that his `algebraic works [\dots]\ were concerned with a problem that had \emph{formally} the same structure as the problem of freedom from contradiction of the classical calculus of logic';\footnote{Letter to \german{Carl Friedrich Gethmann}, \citealp[see][p.~76]{gethmann91}.} compare also his letter to Eckart Menzler-Trott (see \citealp[p.~260]{menzlertrott01}). 

In his letter dated 13 March 1944, he announces: `Subsequently to an algebraic investigation of orthocomplemented semilattices, I am now trying to get out the connection of these questions with the freedom from contradiction of classical logic. [\dots]\ actually I am much more interested into the algebraic side of proof theory than into the purely logical'.\footnote{PL~1-1-131, edited in \citealt[§~6]{neuwirthkonstanz}.} The concept of `orthocomplementation'\footnote{The terminology might be adapted from \citealt{stone36}, where it has a Hilbert space background; today one says `pseudocomplementation'.} (see p.~\pageref{def:orthocomplemented} for its definition) must have been motivated by logical negation from the beginning. On the one hand, such lattices correspond to the calculus of sequents considered by \citet[section~IV]{gentzen36}, who shows that a given derivation can be transformed into a derivation `in which the connectives~$\lor$, $\exists$~and~$\supset$ no longer occur' and provides a proof of consistency for this calculus (see section~\ref{sec:comp-with-gentz} below). On the other hand, note that Lorenzen reviews \citealt{ogasawara39} for the \emph{Zentralblatt}. 

\section{The 1944 manuscript}\label{sec:1944-manuscript}

The result of this investigation can be found in the manuscript `\german{Ein halb\-ordnungstheoretischer Widerspruchsfreiheitsbeweis}'.\footnote{`A proof of freedom from contradiction within the theory of partial order', \german{Oskar-Becker-Nachlass, Philosophisches Archiv, Universität Konstanz}, OB~5-3b-5, \url{https://archive.org/details/lorenzen-ein_halbordnungstheoretischer_widerspruchsfreiheitsbeweis}. The file OB~5-3b consists of documents related to Lorenzen, the oldest being the 1944 manuscript and the youngest a letter from 1951. Lorenzen and Becker are both at Bonn from 1945 to 1956 and have been in close contact since at least 1947: see Lorenzen's letter to Gethmann (in \citealp[p.~77]{gethmann91}).}

We believe that it is the one that he assertedly sends to Wilhelm Ackermann,
Gentzen, Hans Hermes and Heinrich Scholz between March and May 1944, and for which he gets a
dissuasive answer from Gentzen, dated 12 September 1944: `I have looked through your attempt at a consistency proof, not in detail, for which I lack the time. However I say this much: the
consistency of number theory cannot be proven so
simply'.\footnote{The letter is reproduced in
  \citealt[p.~372]{menzlertrott01}, and translated in \citealt{menzlertrott07}.}

Our identification of the manuscript is made on the basis
of the following dating: Lorenzen mentions such a manuscript and its recipients in his
letters to Scholz dated 13 May 1944 and 2 June
1944,\footnote{\german{Heinrich-Scholz-Archiv} and PL~1-1-138.} and in a postcard to Hasse dated 25 July
1945;\footnote{Cod.\ Ms.~H.~Hasse 1:1022, edited in \citealt[§~7]{neuwirthkonstanz}.} a letter by Ackermann
dated 11 November 1946 states that he lost a manuscript by Lorenzen `at the partial destruction of his flat by bombs'.\footnote{`\german{So ist auch ein Manuskript, das Sie mir seiner Zeit zuschickten, bei der teilweisen Zerstörung meiner Wohnung durch Bomben verschwunden}' (PL~1-1-125).} Our identification is also consistent with the content of Lorenzen's letter to Menzler-Trott mentioned above. On the other hand, we have not found any hint at another manuscript by Lorenzen for which it could have been mistaken.\footnote{The `unpublished' manuscript `Ein finiter Logikkalkül' mentioned by \citet[p.~20]{lorenzen48} may be dated to 1947 even if we have not spotted a copy of it: the review given there shows that it corresponds to a thread of research described in a letter to Bernays dated 21 February 1947 (ETH-Bibliothek, Hochschularchiv, Hs 975:2950).} The generalisation of his proof of consistency to ramified type theory is first mentioned in a letter from Scholz to Bernays dated 11 December 1945:\footnote{Hs~975:4111.} it corresponds to the manuscript `Die Widerspruchsfreiheit der klassischen Logik mit verzweigter Typentheorie' and is the future part~II of his 1951 article.

This manuscript renews the relationship between logic and lattice theory: whereas boolean algebras were originally conceived for modeling the classical calculus of propositions, and \foreignlanguage{dutch}{Heyting} algebras for modeling the intuitionistic one, here logic comes at the rescue of lattice theory for studying countably complete pseudocomplemented semilattices.

We have found only three contemporaneous occurrences of the notion of countably complete lattice other than $\sigma$-fields of subsets of a given set used in measure theory: \citet[p.~795; article reviewed by Lorenzen for the \emph{\german{Zentralblatt}}]{birkhoff38} speaks of `$\sigma$-lattice, by analogy with the usual notions of $\sigma$-rings and $\sigma$-fields of sets'; in the appendix \emph{The somen as elements of partially ordered sets} of the posthumously published book \citealt{caratheodory56}; the `$\aleph_1$-lattice' in \citealt{vonneumann37}.

Lorenzen describes a countably complete pseudocomplemented semilattice as a deductive calculus on its own, without any reference to a larger formal framework:\footnote{In contradistinction to the `consequence relation' of \citealt{tarski30} which presupposes set theory.} this conception dates back to the `system of sentences' of \cite{hertz22,hertz23}. The rules of the calculus construct the free countably complete pseudocomplemented semilattice over a given preordered set by taking as axioms the inequalities in the set, by defining inductively formal meets and formal negations, and by introducing inequalities between the formal elements. The introduction rule for formal countable meets, stating that
\[
\text{if $c\lt a_1,c\lt a_2,\dots$, then $c\lt\bigland M$, where $M=(a_1,a_2,\dots)$}
\]
(rule~$c$ on p.~\pageref{rule-c}), stands out: it has an infinity of premisses, so that it is an `$\omega$-rule' in today's terminology. Lorenzen's boldness is most probably due to his training in algebra, where such a rule is very natural, so that when he arrives at a clear constructive understanding of ideal theory, he has also got a clear constructive understanding of the $\omega$-rule.

In ideal theory, \citet[\S\ 4]{MR0033826} defines a system of ideals for a preordered set as the free semilattice generated by it: it consists in the formal meets $a_1\land\dots\land a_m$ of finitely many elements $a_1,\dots,a_m$; this formal element is introduced with the following rules: if $c\lt a_1,\dots,c\lt a_m$, then $c\lt a_1\land\dots\land a_m$; $a_1\land\dots\land a_m\lt a_1,\allowbreak\dots,a_1\land\dots\land a_m\lt a_m$. The $\omega$-rule is the infinitary counterpart of the first rule, and the infinitary counterpart of the second rule is the admissible rule~$\varepsilon$ on p.~\pageref{rule-beta}.

Lorenzen's presentation of elementary number theory can be compared to that of \citealt{goedel33} as follows.
\begin{itemize}
\item Lorenzen starts with `prime formulas', i.e.\ the numerical
  propositions as e.g.\ $1=1''$ or $1+1=1'$. These are preordered by
  material implication and may be combined into compound
  formulas. Lorenzen works in a constructive metatheory, in which infinitely
  many propositions may be supervised if given by a construction,
  e.g.\ the propositions
  $\mathfrak C\rightarrow\mathfrak{A}(1),\mathfrak
  C\rightarrow\mathfrak{A}(2),\dots$, and rule~$c$ on
  p.~\pageref{rule-c-zt} is the rule of introduction of the universal
  quantifier that one may infer from these
  $\mathfrak C\rightarrow(\mathfrak x)\,\mathfrak{A}(\mathfrak x)$.
\item Gödel starts with `elementary formulas', which may also contain
  variables. He works in a finitary metatheory in which only finitely
  many propositions may be supervised, and formalises elementary number
  theory with the universal quantifier handled in a way that is
  equivalent to its usual introduction and elimination rules. Here one
  may construct as in \citealt{goedel31} a predicate
  $\mathfrak{A}(\mathfrak x)$ such that each of the propositions
  $\mathfrak C\rightarrow\mathfrak{A}(1),\mathfrak
  C\rightarrow\mathfrak{A}(2),\dots$ holds, but
  $\mathfrak C\rightarrow(\mathfrak x)\,\mathfrak{A}(\mathfrak x)$
  does not.
\end{itemize}

In elementary number theory, the rule of complete induction plays a central rôle. The statement of this rule is complex from a logical point of view because of the presence of a free variable, of a universal quantifier, or of an implication. The $\omega$-rule appears as an analysis of this complexity: the rule of complete induction is the derivation of $\mathfrak A(1)\rightarrow(\mathfrak x)\,\mathfrak A(\mathfrak x)$ from $\mathfrak A(\mathfrak a)\rightarrow\mathfrak A(\mathfrak a+1)$ with a free variable~$\mathfrak a$; in the latter, replacement of~$\mathfrak a$ by $1,2,\dots$ and the cut rule yield $\mathfrak A(1)\rightarrow\mathfrak A(2),\allowbreak\mathfrak A(1)\rightarrow\mathfrak A(3),\dots$; therefore this rule is a combination of the admissible cut rule~$k$ on p.~\pageref{rule-k} with the $\omega$-rule that derives $\mathfrak A(1)\rightarrow(\mathfrak x)\,\mathfrak A(\mathfrak x)$ from $\mathfrak A(1)\rightarrow\mathfrak A(1),\mathfrak A(1)\rightarrow\mathfrak A(2),\mathfrak A(1)\rightarrow\mathfrak A(3),\dots$\,. Conversely, the only expected uses of the $\omega$-rule correspond to the rule of complete induction and to the rule of introduction of the universal quantifier. The $\omega$-rule has a very simple structure: its premisses are stated without further need of free variables and quantifiers; however, there are infinitely many. Its main feature is that it allows for derivations without detour. 

\citet{sundholm83} and \citet{feferman86} provide a historical account of such rules. \citet{hilbert31,hilbert31b} states an $\omega$-rule with the motivation of, respectively, proving the completeness of arithmetic and the law of excluded middle.\footnote{\cite{hilbert31} states a restricted $\omega$-rule, in the sense that its premisses must be decidable (i.e.\ numerical); he states the axiom of complete induction separately. This is noted in the letter that Bernays addresses to Gödel on 18 January 1931 \citep[pp.~80--91]{goedel03}, where he formulates its unrestricted counterpart. See also Gödel's answer dated 2 April 1931. The $\omega$-rule in \citealt{hilbert31b} is not restricted. Compare \citealt*[pp.~788--805, 967--973, 983--984]{ewaldsiegetal13}.} He declares that it is a `finitary deduction rule', that it has a `rigorously finitary character'. Lorenzen makes no reference to these articles, but, in the 1945 manuscript `Die Widerspruchsfreiheit der klassischen Logik mit verzweigter Typentheorie', he expands on the finitary character of its usage: `One has to persuade oneself at each appearance of this rule that its application occurs to the effect of a ``finitary deduction'', because the proof of freedom from contradiction would otherwise become meaningless'.\footnote{`\german{Man hat sich bei jedem Vorkommen dieser Regeln zu überzeugen, daß ihre Anwendung im Sinne des ``finiten Schließens'' geschieht, weil sonst der Wf-Beweis sinnlos würde}' (`The freedom from contradiction of classical logic with ramified type theory'; a version of this manuscript can be found in \german{Niedersächsische Staats- und Universitätsbibliothek Göttingen}, Cod.\ Ms.~G.~Köthe M~10).}
E.g.\ in the derivation of the rule of complete induction on p.~\pageref{complete-induction}, the infinitely many premisses $\mathfrak A(1)\rightarrow\mathfrak A(1),\mathfrak A(1)\rightarrow\mathfrak A(2),\dots$ must result from a construction whose explanation is finitary, but whose realisation is endless: `For every number~$m$ follows therefrom at once $\mathfrak{A}(1)\rightarrow\mathfrak{A}(m)$ by $m$-fold application of the rule of inference~$k$'. Lorenzen shares this intuitionistic framework with \citet[p.~526]{gentzen36}: `After all, we need not associate the idea of a closed infinite number of individual propositions with this~[$(\mathfrak x)\,\mathfrak A(\mathfrak x)$, where $\mathfrak A$ shall not yet contain an universal or existential quantifier], but can, rather, interpret its sense ``finitistically'' as follows: ``If, starting with~1, we substitute for~$\mathfrak x$ successive natural numbers then, however far we may progress in the formation of numbers, a true proposition results in each case''{}'.

In his letter to Bernays dated 2 April 1931, Gödel points out that such a rule presupposes a framework in which this infinity of premisses may be asserted: `the very complicated and problematical concept ``finitary proof'' is assumed [\dots] without having been made mathematically precise' (see \citealt[p.~97]{goedel03}). This framework is thus an informal one; and, as the proof of consistency rests on its reliability, this framework is to be the intuitionistic one, as \citet[`groupe~D', p.~5]{herbrand31} and \citet[p.~231]{novikoff43} state, i.e.\ the constructive one (\citealt[p.~82]{lorenzen51}). In this sense, a calculus including the $\omega$-rule is of a different nature than a mechanical calculus, where we can check by a finitary process the correctness of a given derivation. In fact, neither the Hilbert program nor Lorenzen's proof of consistency take place in a mechanical formal system, i.e.\ in a system whose objects are finitary and whose derivations are finitary and decidable.

The proof that the calculus thus defined is a countably complete pseudocomplemented semilattice illustrates, as Lorenzen realises a posteriori,\footnote{This is how we interpret the beginning of the second paragraph on p.~\pageref{beginning-2nd}: `Without knowledge of [\dots]'.} that the strategy of Gentzen's dissertation (\textit{\citeyear{gentzen34}}, IV, §~3) for proving the consistency of elementary number theory without complete induction may be maintained for proving the consistency of all of elementary number theory: the introduction rules (rules~$a$ to~$f$ on p.~\pageref{rule-c}) introduce inequalities for formal elements of increasing complexity, i.e.\ no inequality can result from a detour; then the corresponding elimination rules (rules~$\gamma$ to~$\varepsilon$ on p.~\pageref{rule-beta}) are shown to hold by an induction on the complexity of the introduced inequality (in Lorenzen's later terminology, one would say that these rules are shown to be `admissible' and can be considered as resulting from an `inversion principle'); at last transitivity of the preorder, i.e.\ the cut rule (rule~$\beta$ on p.~\pageref{rule-beta}: if $a\leqslant b$ and $b\leqslant c$, then $a\leqslant c$), is established by proving a stronger rule through an induction on the complexity of the cut element~$b$ nested with inductions on the complexity of the derivation of the rule's premisses.

The inductions used here are the ones accurately described by Jacques \citet[pp.~4--5]{herbrand30} after having been emphasised by David \citet[p.~76]{hilbert28}: the first proceeds along the construction of formulas starting from prime formulas through rules, and has no special name (it will be called `formula induction' in \citealt{lorenzen51}); the second proceeds along the construction of theorems starting from prime theorems through deduction rules, and is called `premiss induction'.\footnote{See \citealt{lorenzen39b} for his interest in the foundation of inductive definitions.} 

In other words, Lorenzen starts with a preordered set~$\mathrm{P}$, constructs the free countably complete pseudocomplemented semilattice~$\mathrm{K}$ over~$\mathrm{P}$ and emphasises conservativity, i.e.\ that no more inequalities come to hold among elements of~$\mathrm{P}$ viewed as a subset of~$\mathrm{K}$ than the ones that have been holding before:\footnote{This is exactly the approach of \citet[§~2]{skolem21} for constructing the free lattice over a preordered set, in the course of studying the decision problem for lattices.} one says that $\mathrm{P}$~is embedded into~$\mathrm{K}$ and that the preorder of~$\mathrm{P}$ is embedded into the countably complete preorder of~$\mathrm{K}$.

Then the consistency of elementary number theory with complete induction is established in §~3 by constructing the free countably complete pseudocomplemented semilattice over its `prime formulas', i.e.\ the numerical formulas, viewed as a set preordered by material implication.

Note the presence of rule~$g$ on pp.~\pageref{rule-c} and~\pageref{rule-c-zt}, a contraction rule. This should be put in relation
\begin{itemize}
\item with the rôle of contraction, especially for steps~13.\,5\,1--13.\,5\,3, in Gentzen's proofs of consistency (\citeyear{gentzen36,gentzen74});
\item with the calculus of P. S.~\citet[lemma~6]{novikoff43}, in which contraction may be proved.
\end{itemize}

\section{Comparison with Gentzen's proof of consistency}
\label{sec:comp-with-gentz}

There are similarities and differences with respect to the strategy developed by \german{Gentzen} for proving the consistency of elementary number theory with complete induction. In his first proof, submitted in August 1935, withdrawn and finally published posthumously by Bernays in \citeyear{gentzen74} (after its translation by \citealt{gen1969}), Gentzen defines a concept of reduction procedure for a sequent and shows that such a procedure may be specified for every derivable sequent but not for the contradictory sequent $\rightarrow 1=2$. Let us emphasise two aspects of this concept.
\begin{itemize}
\item If the succedent of the sequent has the form $\forall x\,F(x)$, the following step of the reduction procedure consists in replacing it by $F(n)$, where $n$ is a number to be chosen freely.
\item A reduction procedure is defined as the specification of a sequence of steps for all possible free choices, with the requirement that the reduction terminates for every such choice.
\end{itemize}
In his letter to Bernays dated 4 November 1935,\footnote{Hs~975:1652, translated by \citet[pp.~241--244]{plato17}.} Gentzen visualises a reduction procedure as a tree whose every branch terminates.

The proof that a reduction procedure may be specified for every derivable sequent is by theorem induction. For this, a lemma is needed, claiming that if reduction procedures are known for two sequents $\Gamma\rightarrow D$ and $D,\Delta\rightarrow C$, then a reduction procedure may be specified for their cut sequent $\Gamma,\Delta\rightarrow C$. The proof goes by induction on the construction of the cut formula~$D$ and traces the claim back to the same claim with the same cut formula, but with the sequent $D,\Delta\rightarrow C$ replaced by a sequent $D,\Delta^*\rightarrow C^*$ resulting from it after one or more reduction steps and the cut sequent replaced by $\Gamma,\Delta^*\rightarrow C^*$. By definition of the reduction procedure, this tracing back must terminate eventually.

This last kind of argument may be considered as an infinite descent in the reduction procedure. In his letter to Bernays, Gentzen seems to indicate that this infinite descent justifies an induction on the reduction procedure; as analysed by William W.~\citet{tait15}, this would be an instance of the Bar theorem. But in his following letter, dated 11 December 1935,\footnote{Hs~975:1653, translated by \citet[p.~244]{plato17}.} he writes that `[his] proof is not satisfactory' and announces another proof, to be submitted in February 1936: in it, he defines the concept of reduction procedure for a derivation (and not for a sequent), associates inductively an ordinal to every derivation, and shows that a reduction procedure may be specified for every derivation by an induction on the ordinal.

Let us compare this strategy with Lorenzen's.
\begin{itemize}
\item The free choice is subsumed in a deduction rule, an $\omega$-rule as described above (rules~$c$ and~$j$ on p.~\pageref{rule-c-zt}).\footnote{Compare Bernays' suggestion in his letter to Gentzen dated 9 May 1938, Hs~975:1661, translated by \citet[pp.~254--255]{plato17}.}
\item Elementary number theory is constructed as the cut-free derivations starting from the numerical formulas, so that it is trivially consistent, and the cut rule (rule~$k$ on p.~\pageref{rule-k}) is shown to be admissible: if derivations are known for two sequents $\mathfrak A\rightarrow\mathfrak B$ and $\mathfrak B\rightarrow\mathfrak C$, then a derivation may be specified for their cut sequent $\mathfrak A\rightarrow\mathfrak C$ by a formula induction on the cut formula~$\mathfrak B$ nested with several instances of a theorem induction.
\end{itemize}
In this way, Lorenzen's strategy may be used to realise the endeavour expressed by
\citet{tait15}: `the gap in Gentzen’s argument is filled, not by the
Bar Theorem, but by taking as the basic notion that of a [cut-free]
deduction tree in the first place rather than that of a reduction
tree'. His 1944 proof can thus be seen as a formal improvement on
Gentzen’s 1935 argument, which is all the more remarkable given
Gentzen’s reaction to Lorenzen’s proof.

\section{Comparison with Novikov's proof of consistency}
\label{sec:comp-with-novikov}

\citet{novikoff43} introduces an intuitionistic calculus that contains an $\omega$-rule (rule~6 on p.~233). He defines in §~4 the concept of `regular formula' that expresses that the formula has a cut-free proof, and shows in §~8 that it is an explanation of classical truth. In fact, he proves essentially that cut (`the rule of inference') is admissible. This proof does not use any induction on the cut formula, contrary to Gentzen's and Lorenzen's proofs (see  \citealt{mints91} and \citealt{tupailo92}). In his introduction, Novikov writes: `As a basis, the consistency of which is assumed, the intuitionistic mathematics is taken. From such a point of view it appears to be possible to prove the consistency of [elementary number theory]'.

\section{Mathematical comments}

On p.~\pageref{reductio}, the premiss induction that establishes rule~$\gamma$ is given the form of a reductio ad absurdum, but the reasoning may easily be unraveled into a direct form.

The calculus~$\rmN$ presented on p.~\pageref{N} is in fact common to intuitionistic and classical arithmetic: recall that `the connectives~$\lor$, $\exists$~and~$\supset$ no longer occur'. It may be criticised for its sloppy way of treating variables. 

Furthermore, the introduction of free variables seems dispensable in the presence of an $\omega$-rule. Rule~$j$ and the corresponding elimination rule~$p$ may be omitted from the calculus at the affordable price of giving complete induction the less elegant form $\mathfrak A(1)\Land(\mathfrak x)\,\overline{\mathfrak A(\mathfrak x)\Land\overline{\mathfrak A(\mathfrak x')}}\rightarrow(\mathfrak x)\,\mathfrak A(\mathfrak x)$ as in \citealt{lorenzen62}.

\section{Conclusion}\label{sec:conclusion}

Proof theory continues to focus on measures of complexity by ordinal numbers. The fact that Lorenzen does not resort to ordinals in his proof of consistency should be considered as a feature of his approach.

Lorenzen's article is remarkable for its metamathematical standpoint. A mathematical object is presented as a construction described by rules. A claim on the object is established by an induction that expresses the very meaning of the construction.

The relations between these objects, of the form of an inequality or of an implication, also admit such a presentation: it has the feature that the construction of a relation proceeds as accumulatively (`without detour', i.e.\ cut) as the construction of the formulas appearing in the relation. It is only in a second place that the corresponding elimination rules and the cut rule are shown to be admissible.

In elementary number theory and for the free countably complete pseudocomplemented semilattice, the construction of a relation uses an $\omega$-rule that is stronger than the rule of complete induction but requires infinitely many premisses, so that a relation corresponds to a well-founded tree.

Lorenzen's standpoint holds equally well for a logical calculus and for a lattice: `logical calculuses\footnote{We prefer this plural with \citet{curry58}.} are semilattices or lattices' (\citealp[p.~89]{lorenzen51}). The consistency of the logical calculus of elementary number theory is recognised as a consequence of the following fact: a preordered set embeds into the free countably complete pseudocomplemented semilattice generated by it in a conservative way.

Other reflections on the philosophical significance of Lorenzen's approach to logic are addressed by Matthias \citet{wille13,wille16}.

\section*{Acknowledgments}

We thank Brigitte Parakenings for having provided ideal working conditions and her expertise at \german{Philosophisches Archiv} of \german{Universität Konstanz}, and Henri Lombardi and Jan von Plato for helpful discussions. This research has been supported through the program `Research in pairs' of Mathematisches Forschungsinstitut Oberwolfach in 2016 and through the hospitality of the university of Gothenburg.

\clearpage\appendix

% \setlength\abovedisplayskip{5pt}\setlength\belowdisplayskip{5pt}\setlength\parskip{0pt}\setlength\medskipamount{4pt}
% \the\abovedisplayskip\the\parskip\the\medskipamount\the\topsep\the\itemsep\the\partopsep\the\parsep
\renewcommand\bigland{\mathop{\textstyle\bigwedge}\limits}
\begin{Parallel}[p]{}{}
\ParallelLText{%
\rightline{[P. LORENZEN]\qquad\qquad}\medskip
  
\begin{center}
\uline{Ein halbordnungstheoretischer Widerspruchsfreiheitsbeweis.}
\end{center}

Die Dissertation von G. Gentzen enthält einen Wf-beweis
der reinen Zahlentheorie ohne vollständige Induktion, der
auf dem folgenden Grundgedanken beruht: jede herleitbare
Sequenz muß sich auch ohne Umwege herleiten lassen, sodaß
während der Herleitung nur die Verknüpfungen eingeführt werden,
die unbedingt notwendig sind, nämlich diejenigen, die in der
Sequenz selbst enthalten sind. In dem Wf-beweis der Zahlen\bs
theorie mit vollständiger Induktion tritt dieser Grundgedanke
gegenüber anderen zurück. Ich möchte jedoch im folgenden
zeigen, daß er allein genügt, auch diese Wf.\ zu erhalten.

Ohne Kenntnis der Dissertation von Gentzen bin ich auf
diese Möglichkeit auf Grund einer halbordnungstheoretischen
Frage gekommen. Diese lautete: wie läßt sich eine halb\bs
geordnete Menge in einen orthokomplementären vollständigen
Halbverband einbetten? Im allgemeinen sind mehrere solche
Einbettungen möglich --~unter den möglichen Einbettungen ist
aber eine ausgezeichnet, nämlich die, welche sich in jede
andere homomorph abbilden läßt. Die Existenz dieser ausge\bs
zeichneten Einbettung wird in §~2 bewiesen.

Um hieraus in §~3 den gesuchten Wf-beweis zu erhalten,
ist nur noch eine Übersetzung des halbordnungstheoretischen
Beweises in die logistische Sprache notwendig. Denn der Kalkül,
den wir betrachten und auf den sich die üblichen Kalküle zu\bs
rückführen lassen, ist in der ausgezeichneten Einbettung der
halbgeordneten Menge der zahlentheoretischen Primformeln ent\bs
halten.\bruch{2}\medskip

\uline{§\enspace1.} Eine Menge~$\rmM$ heißt \red{halbgeordnet}, wenn in~$\rmM$ eine zwei\bs
stellige Relation~$\leqslant$ definiert ist, sodaß für die Elemente
$a,b,\dots$ von~$\rmM$ gilt:%
\[\begin{gathered}%
  a\leqslant a\\%
  a\leqslant b,\;b\leqslant c\quad\Rightarrow\quad a\leqslant c\text.%
\end{gathered}\]
Gilt $a\leqslant b$ und $b\leqslant a$, so schreiben wir $a\equiv b$.

Gilt $a\leqslant x$ für jedes~$x\in \rmM$, so schreiben wir $a\leqslant{}$. Ebenso
schreiben wir ${}\leqslant a$, wenn $x\leqslant a$ für jedes~$x$ gilt. (${}\leqslant{}$ bedeutet
also, daß $x\leqslant y$ für jedes~$x,y\in \rmM$ gilt.)

Eine halbgeordnete Menge~$\rmM$ heißt \red{Halbverband}, wenn es zu jedem
$a,b\in \rmM$ ein~$c\in \rmM$ gibt, sodaß für jedes~$x\in \rmM$ gilt%
\[x\leqslant a,\;x\leqslant b\quad\Longleftrightarrow\quad x\leqslant c\text.\]%
$c$ heißt die \emph{Konjunktion} von~$a$ und~$b$: $c\equiv a\land b$.

Ein Halbverband~$\rmM$ heißt \red{orthokomplementär}, wenn es zu jedem~$a\in \rmM$
ein~$b\in \rmM$ gibt, so daß für jedes~$x\in \rmM$ gilt%
\[a\land x\leqslant{}\quad\Longleftrightarrow\quad x\leqslant b\text.\]%
$b$~heißt das \emph{Orthokomplement} von~$a$: \ $b\equiv\neg{a}$.

Ein Halbverband~$\rmM$ heißt \red{$\omega$-vollständig}, wenn es zu jeder abzähl\bs
baren Folge~$\itM=a_1,a_2,\dots$ in~$\rmM$ ein~$c\in \rmM$ gibt, so daß für jedes\nbsp
$x\in \rmM$ gilt:%
\[(\text{für jedes~$n$: }x\leqslant a_n)\quad\Longleftrightarrow\quad x\leqslant c\text.\]%
$c$~heißt die \emph{Konjunktion der Elemente von}~$\itM$: \ $c\equiv\bigland_na_n\equiv\bigland_{\itM}{}$.

Sind~$\rmM$ und~$\rmM'$ halbgeordnete Mengen, so heißt~$\rmM$ \red{ein Teil} von $\rmM'$,
wenn~$\rmM$ Untermenge von~$\rmM'$ ist und für jedes~$a,b\in \rmM$ genau dann
$a\leqslant b$ in~$\rmM'$ gilt, wenn $a\leqslant b$ in~$\rmM$ gilt.

Sind~$\rmM$ und~$\rmM'$ halbgeordnete Mengen, so verstehen wir unter
einer \red{Abbildung} von~$\rmM$ in~$\rmM'$ eine Zuordnung, die jedem~$a\in \rmM$
ein~$a'\in \rmM'$ zuordnet, so daß gilt%
\[a\equiv b\quad\Rightarrow\quad a'\equiv b'\text.\bruch{3}\]

Sind~$\rmM$ und~$\rmM'$ orthokomplementäre $\omega$-vollständige Halb\bs
verbände, so verstehen wir unter einem \red{Homomorphismus} von~$\rmM$
in~$\rmM'$ eine Abbildung~$\to$ von~$\rmM$ in~$\rmM'$, so daß für jedes~$a,b\in \rmM$
und~$a',b'\in \rmM'$ mit $a\to a'$ und $b\to b'$ gilt:%
\[%
  \begin{aligned}
    a\land b&\to a'\land b'\\%
    \neg a&\to\neg{a'}\text.%
  \end{aligned}
\]%
Ferner soll für jede Folge~$\itM=a_1,a_2,\dots$ in~$\rmM$ und
$\itM'=a'_1,a'_2,\dots$ in~$\rmM'$ mit $a_n\to a'_n$ gelten:%
\[\bigland_{\itM}{}\to\bigland_{\itM'}{}\text.\]

Wir wollen jetzt beweisen, daß es zu jeder halbgeordneten
Menge~$\rmP$ einen orthokomplementären $\omega$-vollständigen Halbver\bs
band~$\rmK$ gibt, so daß%
\begin{itemize}%
\item[1)]$\rmP$~ein Teil von~$\rmK$ ist,% [d.\ h.\ $\rmK$~enthält~$\rmP$ als Teil.]
\item[2)]$\rmK$ in jeden orthokomplementären $\omega$-vollständi\bs
gen Halbverband, der $\rmP$~als Teil enthält,
homomorph abbildbar ist.%
\end{itemize}%
Wäre $\rmK'$~ein weiterer orthokomplementärer $\omega$-vollständiger
Halbverband, der die Bedingungen 1)~und~2) erfüllt, so gäbe
es eine Zuordnung, durch die $\rmK$~in~$\rmK'$ und $\rmK'$~in~$\rmK$ homomorph
abgebildet würde, d.\ h.\ $\rmK$~und~$\rmK'$ wären \red{isomorph}. $\rmK$~ist also
durch die Bedingungen 1)~und~2) bis auf Isomorphie eindeutig
bestimmt. Wir nennen $\rmK$~\red{den ausgezeichneten orthokomplemen\bs
tären $\omega$-vollständigen Halbverband über~$\rmP$}.\medskip

\uline{§\enspace2.} \uline{Satz}: \emph{Über jeder halbgeordneten Menge gibt es den aus\bs
gezeichneten orthokomplementären\/ $\omega$-vollständigen Halbverband.}

Wir konstruieren zu der halbgeordneten Menge~$\rmP$ eine Menge~$\rmK$
auf folgende Weise:%  
\begin{enumerate}[label=\arabic*)]%
\item $\rmK$~enthalte die Elemente von~$\rmP$. (Diese nennen wir die
\red{Primelemente} von~$\rmK$.)\bruch{4}
\item $\rmK$~enthalte mit endlich vielen Elementen~$a_1,a_2,\dots,a_n$
auch die hieraus gebildete \red{Kombination} als Element.
(Diese bezeichnen wir durch~$a_1\land a_2\land\cdots\land a_n$.)%
\item $\rmK$~enthalte mit jedem Element~$a$ auch ein Element~$\neg a$.%
\item $\rmK$~enthalte mit jeder abzählbaren Folge~$\itM$ auch ein
Element~$\bigland_\itM{}$.
\end{enumerate}

Jedes Element von~$\rmK$ läßt sich also eindeutig als Kombination\nbsp
$a_1\land a_2\land\cdots\land a_n$ von Primelementen und Elementen der Form~$\neg a$
oder~$\bigland_{\itM}{}$ schreiben.

Wir definieren eine Relation~$\leqslant$ in~$\rmK$ auf folgende Weise:%
\begin{itemize}%
\item[1)]Für Primelemente~$p,q$ gelte $p\leqslant q$ in~$\rmK$, wenn~$p\leqslant q$ in~$\rmP$
gilt. (Diese Relationen nennen wir die Grundrelationen.)%
\item[2)]Es soll jede Relation~$\leqslant$ in~$\rmK$ gelten, die sich aus den
Grundrelationen mit Hilfe der folgenden Regeln herlei\bs
ten läßt:%
\[%
  \begin{gathered}%
    \begin{aligned}[t]%
      &{\begin{prooftree}%
          \Hypo{c\leqslant a}\Hypo{c\leqslant b}\Infer[left label=$a)$]2{c\leqslant a\land b}%
        \end{prooftree}}\\%
      &{\begin{prooftree}%
          \Hypo{a\land c\leqslant\hphantom a}\Infer[left label=$b)$]1{\hphantom{a\land{}}c\leqslant\neg a}%
        \end{prooftree}}\\%
      &{\begin{prooftree}[separation=1ex]%
          \Hypo{c\leqslant a_1,}\Hypo{\dots,}\Hypo{c\leqslant a_n,}\Hypo{\dots}\Infer[left label=$c)$]4{c\leqslant\bigland_{\itM}}%
        \end{prooftree}}\\%
    \end{aligned}\quad%
    \begin{aligned}[t]%
      &{\begin{prooftree}%
          \Hypo{\hphantom{b\land{}}a\leqslant c}\Infer[left label=$d{})$]1{a\land b\leqslant c}%
        \end{prooftree}}\\%
      &{\begin{prooftree}%
          \Hypo{\hphantom{\neg b\land{}}a\leqslant b}\Infer[left label=$e)$]1{a\land\neg b\leqslant c}%
        \end{prooftree}}\\%
      &{\begin{prooftree}%
          \Hypo{a_n\land b\leqslant c}\Infer[left label=$f)$]1{\hphantom{a_n}\llap{$\bigland_{\itM}$}\land b\leqslant c}%
        \end{prooftree}}%
    \end{aligned}\\%
    (\itM=a_1,a_2,\dots)\\%
    {\begin{prooftree}%
        \Hypo{a\land a\land b\leqslant c}\Infer[left label=$g)$]1[\bruch{5}]{\hphantom{a\land{}}a\land b\leqslant c}%
      \end{prooftree}}%
  \end{gathered}%
\]%
\end{itemize}%
Wir nennen die Relationen über dem Strich die \red{Prämissen}
der Relation unter dem Strich.

Wir haben jetzt zunächst zu zeigen, daß~$\rmK$ ein orthokomple\bs
mentärer $\omega$-\hspace{0pt}vollständiger Halbverband bezügl.\ der Relation~$\leqslant$
ist. Dazu müssen wir beweisen%
\[%
  \begin{aligned}%
    \alpha)&&a\leqslant a&\\%
    \beta)&&a\leqslant b,\;b\leqslant c\quad&\Rightarrow\quad a\leqslant c\\%
    \gamma)&&c\leqslant a\land b\quad&\Rightarrow\quad c\leqslant a\\%
    \delta)&&c\leqslant\neg a\quad&\Rightarrow\quad a\land c\leqslant{}\\%
    \varepsilon)&&c\leqslant\bigland_{\itM}\quad&\Rightarrow\quad c\leqslant a_n&&(\itM=a_1,a_2,\dots)%
  \end{aligned}%
\]%
Diese Eigenschaften zusammen mit~$a)$, $b)$~und~$c)$ drücken nämlich
aus, daß~$\rmK$ ein orthokomplementärer $\omega$-\hspace{0pt}vollständiger Halb\bs
verband ist.

$\alpha)$ gilt für Primelemente. Gilt $\alpha)$ für~$a$ und~$b$, so auch für\nbsp
$a\land b$ wegen%
\begin{prooftree*}%
  \Hypo{\hphantom{b\land{}}a\leqslant a}\Infer1{a\land b\leqslant a}\Hypo{\hphantom{a\land{}}b\leqslant b}\Infer1{a\land b\leqslant b}\Infer2{a\land b\leqslant a\land b}%
\end{prooftree*}
Gilt $\alpha)$ für jedes~$a_n\in\itM$, so auch für~$\bigland_{\itM}$ wegen%
\begin{prooftree*}%
  \Hypo{a_1\leqslant a_1}\Infer1{{\bigland_{\itM}}\leqslant a_1}\Hypo{\cdots}\Infer[no rule]1{\vrule height5pt width0pt}%\rewrite{\raise2ex\box\treebox}
  \Hypo{a_n\leqslant a_n}\Infer1{{\bigland_{\itM}}\leqslant a_n}\Hypo{\cdots}\Infer[no rule]1{\vrule height5pt width0pt}%\rewrite{\raise2ex\box\treebox}
  \Infer4{{\bigland_{\itM}}\leqslant\bigland_{\itM}}%
\end{prooftree*}
Gilt $\alpha)$ für~$a$, so auch für~$\neg a$, wegen%
\begin{prooftree*}%
  \Hypo{\hphantom{\neg a\land{}}a\leqslant a}\Infer1{a\land\neg a\leqslant\hphantom{a}}\Infer1{\hphantom{a\land{}}\neg a\leqslant\neg a}%
\end{prooftree*}%
Dadurch ist $\alpha)$ allgemein bewiesen.\bruch{6}

Da $\beta)$ am schwierigsten zu beweisen ist, nehmen wir zunächst $\gamma)$.

Um $\gamma)$ zu beweisen, haben wir zu zeigen, daß, wenn~$c\leqslant a\land b$
herleitbar ist, dann auch stets $c\leqslant a$ herleitbar sein muß.

% \looseness=-1
Wir führen den Beweis indirekt durch eine \emph{transfinite Induk\bs
tion}. Es sei~$c\leqslant a\land b$ herleitbar, aber nicht~$c\leqslant a$. Der letzte
Schritt der Herleitung von $c\leqslant a\land b$ kann dann nicht sein
${\begin{prooftree}[center]\Hypo{c\leqslant a}\Hypo{c\leqslant b}\Infer2{c\leqslant a\land b}\end{prooftree}}$ ebenfalls nicht ${\begin{prooftree}[center]\Hypo{\hphantom{\neg{c_2}\land{}}c_1\leqslant\rlap{$c_2$}\hphantom{a\land b}}\Infer1{c_1\land\neg{c_2}\leqslant a\land b}\end{prooftree}}$ ($c=c_1\land\neg{c_2}$)
da dann sofort ${\begin{prooftree}[center]\Hypo{\hphantom{\neg{c_2}\land{}}c_1\leqslant c_2}\Infer1{c_1\land\neg{c_2}\leqslant\rlap{$a$}\hphantom{c_2}}\end{prooftree}}$ herleitbar wäre.

Für den letzten Schritt bleiben nur die Möglichkeiten%
\[%
  \begin{gathered}%
    {\begin{prooftree}%
        \Hypo{\hphantom{c_2\land{}}c_1\leqslant a\land b}\Infer1{c_1\land c_2\leqslant a\land b}%
      \end{prooftree}}\qquad\qquad%
      {\begin{prooftree}%
        \Hypo{c_1\land c_1\land c_2\leqslant a\land b}\Infer1[$(c=c_1\land c_2)$]{\hphantom{c_1\land{}}c_1\land c_2\leqslant a\land b}%
      \end{prooftree}}\\%
      {\begin{prooftree}%
          \Hypo{c_1\land c'\leqslant a\land b}\Hypo{\cdots}\Hypo{c_n\land c'\leqslant a\land b}\Hypo{\cdots}\Infer4[$\left(%
        \begin{aligned}%
          \itM&=a_1,a_2,\dots\\%
          c&=\smash{\bigland_{\itM}}\land c'%
        \end{aligned}%
\right)$]{{\bigland_{\itM}}\land c'\leqslant a\land b}%
        \end{prooftree}}%
  \end{gathered}%
\]

Hier muß jetzt~$c_1\leqslant a$ bezw.\ $c_1\land c_1\land c_2\leqslant a$ bezw.\ für min\bs
destens ein~$n$ \ $c_n\land c'\leqslant a$ nicht herleitbar sein, da sonst sofort
$c\leqslant a$ herleitbar wäre. In der Herleitung von~$c\leqslant a\land b$ wäre also
schon für eine Prämisse die Behauptung~$\gamma)$ falsch. Gehe ich in
der Herleitung von einer Relation zu einer Prämisse über, von
dieser wieder zu einer Prämisse usw., so bin ich nach endlich
vielen Schritten bei einer Grundrelation. Wir erhielten also
eine Grundrelation, für die die Behauptung~$\gamma)$ falsch wäre.
Da dieses aber unmöglich ist, ist damit $\gamma)$~bewiesen.

% \looseness=-1
Wir nennen die Induktion, die wir hier durchgeführt haben,
eine \red{Prämisseninduktion}.

Mit Hilfe von Prämisseninduktionen verläuft der Beweis für\nbsp
$\delta)$ und~$\varepsilon)$ ebenso einfach wie für~$\gamma)$, so daß ich hierauf
nicht weiter eingehe.\bruch{7}

% \looseness=-1
Es bleibt nur noch $\beta)$ zu zeigen. Statt dessen beweisen wir die
stärkere Behauptung%
\[\zeta)\quad a\leqslant b,\; b\land b\land\cdots\land b\land c\leqslant d\quad \Rightarrow\quad a\land c\leqslant d\]%
um hierauf Prämisseninduktionen anwenden zu können.

Es seien zunächst~$b$, $c$~und~$d$ Primelemente. Dann gilt~$\zeta)$ 
für jede Grundrelation~$a\leqslant b$. Wir nehmen als \emph{Induktionsvor\bs
aussetzung} an, daß $\zeta)$~\emph{für jede Prämisse von}~$a\leqslant b$ gelte.

% \looseness=-1
Da $b$~ein Primelement ist, kann der letzte Schritt der Her\bs
leitung von~$a\leqslant b$ nur sein:%
\[%
  \begin{gathered}%
    {\begin{prooftree}%
        \Hypo{\hphantom{a_2\land{}}a_1\leqslant b}\Infer%[left label=(1)]
        1{a_1\land a_2\leqslant b}%
      \end{prooftree}}\qquad%
    {\begin{prooftree}%
        \Hypo{a_1\land a_1\land a_2\leqslant b}\Infer%[left label=(2)]
        1[$(a=a_1\land a_2)$]{\hphantom{a_1\land{}}a_1\land a_2\leqslant b}%
      \end{prooftree}}\\%
    {\begin{prooftree}%
        \Hypo{\hphantom{a_2\land{}}a_1\leqslant a_2\vphantom{a'}}\Infer%[left label=(3)]
        1[$(a=a_1\land\neg{a_2})$]{a_1\land\neg{a_2}\leqslant \rlap{$b$}\hphantom{a_2}\vphantom{\bigland_{\itM}}}%
      \end{prooftree}}\qquad%
    {\begin{prooftree}%
        \Hypo{a_n\land a'\leqslant b}\Infer%[left label=(4)]
        1[$\left(%
          \begin{aligned}%
            \itM&=a_1,a_2,\dots\\%
            a&=\smash{\bigland_{\itM}}\land a'%
          \end{aligned}\right)$]%
        {{\bigland_{\itM}}\land a'\leqslant b}%
      \end{prooftree}}    
  \end{gathered}%
\]

Nach der Induktionsvoraussetzung ist dann $a_1\land c\leqslant d$ bezw.\bsp
$a_1\land a_1\land a_2\land c\leqslant d$ bezw.\ $a_n\land a'\land c\leqslant d$ herleitbar.
In jedem Falle ist sofort~$a\land c\leqslant d$ herleitbar, ebenso aus $a_1\leqslant a_2$ 
wegen%
\begin{prooftree*}%
  \Hypo{\hphantom{a_2\land{}}a_1\leqslant a_2}\Infer1{a_1\land\neg{a_2}\leqslant\rlap{$d$}\hphantom{a_2}}\Infer1{\hphantom{a_1\land a_2}\llap{$a\land c$}\leqslant\rlap{$d$}\hphantom{a_2}}%
\end{prooftree*}

Damit ist $\zeta)$~bewiesen für Primelemente~$b$, $c$~und~$d$.

% \looseness=-1
Jetzt sei nur noch~$b$ ein Primelement. Dann gilt also~$\zeta)$
für beliebiges~$a$ und Primelemente~$c$,~$d$. Eine \emph{Prämissen\bs
induktion} ergibt jetzt, daß $\zeta)$ für jede Relation $b\land b\land\cdots\land b\land c\leqslant d$
gilt. Jede Prämisse von $b\land b\land\cdots\land b\land c\leqslant d$ hat nämlich wieder
die Form $b\land\cdots\land b\land c\leqslant d$. Damit ist $\zeta)$ allgemein für \emph{Prim\bs
elemente}~$b$ bewiesen.

Gilt $\zeta)$ für Elemente~$b_1$ und~$b_2$, so auch ersichtlich für $b_1\land b_2$.
Gilt $\zeta)$ für jedes~$b_n\in\itM$, so auch für~$b=\bigland_{\itM}$. (Beweis durch
Prämisseninduktion: ${\bigland_{\itM}}\land{\bigland_{\itM}}\land\dots\land{\bigland_{\itM}}\land c\leqslant d$ kann folgende
Prämisse haben: $b_n\land{\bigland_{\itM}}\land\cdots\land{\bigland_{\itM}}\land c\leqslant d$. Nach Induktionsvoraus\bruch{8}%
setzung gilt dann $a\leqslant{\bigland_{\itM}},\;b_n\land{\bigland_{\itM}}\land\cdots\land{\bigland_{\itM}}\land c\leqslant d\quad\Rightarrow\quad b_n\land a\land c\leqslant d$.
Da $\zeta)$ aber auch für $b=b_n$ vorausgesetzt ist, und wegen%
\[a\leqslant{\bigland_{\itM}}\quad\Rightarrow\quad a\leqslant b_n\]%
gilt auch $a\leqslant{\bigland_{\itM}},\;b_n\land a\land c\leqslant d\;\Rightarrow\; a\land a\land c\leqslant d$.
Aus $a\land a\land c\leqslant d$ ist aber $a\land c\leqslant d$ herleitbar. Jede andere
Prämisse von ${\bigland_{\itM}}\land{\bigland_{\itM}}\land\cdots\land{\bigland_{\itM}}\land c\leqslant d$ ist trivial.)

Gilt $\zeta)$ für~$b$, so auch für~$\neg b$. (Beweis durch Prämisseninduk\bs
tion: $\neg b\land\neg b\land\cdots\land\neg b\land c\leqslant d$ kann die folgende Prämisse haben:
$\neg b\land\cdots\land\neg b\land c\leqslant b$. Dann gilt nach Induktionsvoraussetzung%
\[a\leqslant\neg b,\;\neg b\land\cdots\land\neg b\land c\leqslant b\quad\Rightarrow\quad a\land c\leqslant b\text.\]%
Da $\zeta)$ auch für~$b$ vorausgesetzt ist, gilt auch%
\[a\land c\leqslant b,\;a\land b\leqslant d\quad\Rightarrow\quad a\land a\land c\leqslant d\text.\]%
Also gilt auch $a\leqslant\neg b,\;\neg b\land\cdots\land\neg b\land c\leqslant b\quad\Rightarrow\quad a\land c\leqslant d$
wegen $a\leqslant\neg b\quad\Rightarrow\quad a\land b\leqslant d$.
Jede andere Prämisse ist wieder trivial.)

Also ist $\zeta)$ allgemein gültig. Damit ist bewiesen, daß $\rmK$~ein
orthokomplementärer $\omega$-\hspace{0pt}vollständiger Halbverband ist.

$\rmP$~ist ein Teil von~$\rmK$, da%
\[p\leqslant q\text{ in }\rmP\quad\Longleftrightarrow\quad p\leqslant q\text{ in }\rmK\]%
gilt. Wir haben uns dazu zu überzeugen, daß keine Relation
$p\leqslant q$ in~$\rmK$ herleitbar ist, die nicht schon in~$\rmP$ gilt. Das ist
aber selbstverständlich, da keine der Regeln außer~$g)$ über\bs
haupt Relationen~$p\leqslant q$ unter dem Strich liefert. Eine Her\bs
leitung einer Relation~$p\leqslant q$ kann also nur die Regeln~$d{})$ und~$g)$
benutzen. Mit diesen sind aber nur die Grundrelationen her\bs
leitbar.

Zum Beweis unseres Satzes bleibt jetzt noch zu zeigen, daß
sich $\rmK$~in jeden anderen orthokomplementären $\omega$-\hspace{0pt}vollständigen
Halbverband~$\rmK'$, der~$\rmP$ als Teil enthält, homomorph abbilden\bruch{9}
läßt. Diese Abbildung definieren wir durch%
\begin{enumerate}[label=\arabic*)]%
\item für Primelemente~$p$ gilt~$p\to p$,%
\item ferner soll gelten%
  \[%
    \begin{aligned}%
      a\to a',\;b\to b'\quad&\Rightarrow\quad a\land b\to a'\land b'\\%
      a\to a'\quad&\Rightarrow\quad\neg a\to\neg{a'}\\%
      a_n\to a'_n\quad&\Rightarrow\quad\bigland_{\itM}\to\bigland_{\itM'}\quad%
      \left(%
        \begin{aligned}%
          \itM&=a_1,a_2,\dots\\%
          \itM'&=a'_1,a'_2,\dots%
        \end{aligned}%
      \right)%    
    \end{aligned}%
  \]%
\end{enumerate}%
Dadurch wird ersichtlich ein Homomorphismus definiert, denn
es gilt für~$a\to a'$ und~$b\to b'$ stets~$a\leqslant b\quad\Rightarrow\quad a'\leqslant b'$.

Jede Herleitung von $a\leqslant b$ beweist nämlich sofort auch $a'\leqslant b'$,
da die Herleitungsschritte~$a)$~-~$g)$ in jedem orthokomplemen\bs
tären $\omega$-vollständigen Halbverband stets richtig sind.\medskip

\uline{§\enspace3.} Um aus dem im §~2 bewiesenen Satz die Widerspruchsfrei\bs
heit der reinen Zahlentheorie mit vollständiger Induktion
beweisen zu können, benutzen wir die folgende Formalisierung.
Als Primformeln nehmen wir die Zeichen für zahlentheoretische
Prädikate~$\rmA(\dots)$, $\rmB(\dots)$,~\ldots\ mit den Zahlen $1,1',1'',\dots$
als Argumenten, z.\ B. $1 = 1''$, $1+1=1'$.

Diese Primformeln~$\mathfrak P,\mathfrak Q,\dots$ bilden eine halbgeordnete Menge,
wenn wir $\mathfrak P\rightarrow\mathfrak Q$ setzen, falls das Prädikat~$\mathfrak P$ das Prädikat~$\mathfrak Q$
impliziert. Zu den Grundrelationen $\mathfrak P\rightarrow\mathfrak Q$ nehmen wir auch
noch die Relationen der Form ${}\rightarrow\mathfrak P$, $\mathfrak P\rightarrow{}$, ${}\rightarrow{}$ hinzu,
soweit sie inhaltlich richtig sind.

Über dieser halbgeordneten Menge~$\rmP$ der Primformeln konstru\bs
ieren wir jetzt wie in §~2 den ausgezeichneten orthokom\bs
plementären $\omega$-\hspace{0pt}vollständigen Halbverband. Wir benutzen dazu
die logistischen Zeichen, also~$\rightarrow$ statt~$\leqslant$, $\Land$ statt $\land$.

Zu den Formeln gehören also die Primformeln, mit~$\mathfrak A$ und~$\mathfrak B$
auch $\mathfrak A\Land\mathfrak B$, mit~$\mathfrak A$ auch~$\overline{\mathfrak A}$. Die Konjunktion abzählbarer\bruch{10}
Folgen beschränken wir auf die Folgen der Form $\mathfrak A(1),\mathfrak A(1'),\dots$
Diese Konjunktion bezeichnen wir durch $(\mathfrak x)\,\mathfrak A(\mathfrak x)$.

Ferner führen wir noch freie Variable $\mathfrak a$ = $a$, $b$,~\ldots\ ein durch
folgende Schlußregel:%
\[
  \pbox{214pt}{sind $A(1),A(1'),\dots$ herleitbare Relationen, so\\%
    soll auch $A(\mathfrak a)$ herleitbar sein.}
\]

Hierdurch werden die Beweise von §~2 nur unwesentlich modi\bs
fiziert. Wir erhalten insgesamt einen Kalkül~$\rmN$ mit den fol\bs
genden Schlußregeln%
\[%
  \begin{gathered}%
    \begin{aligned}[t]%
      &{\begin{prooftree}%
          \Hypo{\mathfrak C\rightarrow\mathfrak A}\Hypo{\mathfrak C\rightarrow\mathfrak B}\Infer[left label=$a)$]2{\mathfrak C\rightarrow\mathfrak A\Land\mathfrak B}%
        \end{prooftree}}\\%
      &{\begin{prooftree}%
          \Hypo{\mathfrak A\Land \mathfrak C\rightarrow\hphantom{\overline{\mathfrak A}}}\Infer[left label=$b)$]1{\hphantom{\mathfrak A\Land{}}\mathfrak C\rightarrow\overline{\mathfrak A}}%
        \end{prooftree}}\\%
      &{\begin{prooftree}%
          \Hypo{\mathfrak C\rightarrow\mathfrak A(1)}\Hypo{\cdots}\Hypo{\mathfrak C\rightarrow\mathfrak A(n)}\Hypo{\cdots}\Infer[left label=$c)$]4{\mathfrak C\rightarrow(\mathfrak x)\,\mathfrak A(\mathfrak x)}%
        \end{prooftree}}\\%
    \end{aligned}\quad%
    \begin{aligned}[t]%
      &{\begin{prooftree}%
          \Hypo{\hphantom{\mathfrak B\Land{}}\mathfrak A\rightarrow \mathfrak C}\Infer[left label=$d{})$]1{\mathfrak A\Land\mathfrak B\rightarrow \mathfrak C}%
        \end{prooftree}}\\%
      &{\begin{prooftree}%
          \Hypo{\hphantom{\overline{\mathfrak B}\Land{}}\mathfrak A\rightarrow\mathfrak B}\Infer[left label=$e)$]1{\mathfrak A\Land\overline{\mathfrak B}\rightarrow\mathfrak C}%
        \end{prooftree}}\\%
      &{\begin{prooftree}%
          \Hypo{\hphantom{(\mathfrak x)\,\mathfrak A(\mathfrak x)}\llap{$\mathfrak A(n)$}\Land\mathfrak B\rightarrow \mathfrak C}\Infer[left label=$f)$]1{(\mathfrak x)\,\mathfrak A(\mathfrak x)\Land\mathfrak B\rightarrow \mathfrak C}%
        \end{prooftree}}%
    \end{aligned}\\%
    {\begin{prooftree}%
        \Hypo{\mathfrak A\Land\mathfrak A\Land\mathfrak B\rightarrow \mathfrak C}\Infer[left label=$g)$]1{\hphantom{\mathfrak A\Land{}}\mathfrak A\Land\mathfrak B\rightarrow \mathfrak C}%
      \end{prooftree}}\\%
    {\begin{prooftree}%
        \Hypo{\mathfrak A\Land\mathfrak B\rightarrow\mathfrak C}\Infer[left label=$h)$]1{\mathfrak B\Land\mathfrak A\rightarrow\mathfrak C}%
      \end{prooftree}}\qquad%
        {\begin{prooftree}%
        \Hypo{\mathfrak A\Land(\mathfrak B\Land\mathfrak C)\rightarrow\mathfrak D}\Infer[left label=$i)$]1{(\mathfrak A\Land\mathfrak B)\Land\mathfrak C\rightarrow\mathfrak D}%
      \end{prooftree}}\\%
    {\begin{prooftree}%
          \Hypo{A(1)}\Hypo{\cdots}\Hypo{A(n)}\Hypo{\cdots}\Infer[left label=$j)$]4[\bruch{11}]{A(\mathfrak a)}%
        \end{prooftree}}%
  \end{gathered}%
\]%

Die Schlußregeln~$h)$ und~$i)$ waren in §~2 überflüssig, da wir
dort~$a\land b\land c\dots$ sofort als Zeichen für die Kombination von
$a,b,c,\dots$ eingeführt haben.

Der Beweis in §~2 liefert jetzt das folgende Ergebnis:
Der Kalkül~$\rmN$ ist widerspruchsfrei, z.\ B. ist die leere Re\bs
lation ${}\rightarrow{}$ nicht herleitbar, da nur die inhaltlich
richtigen Relationen in~$\rmP$ gelten und $\rmP$~ein Teil von~$\rmN$ ist.
Zu dem Kalkül~$\rmN$ können die folgenden Schlußregeln hinzuge\bs
nommen werden, ohne daß die Menge der herleitbaren Relati\bs
onen vergrößert wird:%
\[%
  \begin{gathered}%
    {\begin{prooftree}%
        \Hypo{\mathfrak A\rightarrow\mathfrak B}\Hypo{\mathfrak B\rightarrow\mathfrak C}\Infer[left label=$k)$]2{\mathfrak A\rightarrow\mathfrak C}%
      \end{prooftree}}\\%
    \begin{aligned}%
      {\begin{prooftree}%
          \Hypo{\mathfrak C\rightarrow\mathfrak A\Land\mathfrak B}\Infer[left label=$l)$]1{\mathfrak C\rightarrow\mathfrak A\hphantom{{}\Land\mathfrak B}}%
        \end{prooftree}}&&\qquad&%
      {\begin{prooftree}%
          \Hypo{\mathfrak C\rightarrow\mathfrak A\Land\mathfrak B}\Infer[left label=$m)$]1{\mathfrak C\rightarrow\mathfrak B\hphantom{{}\Land\mathfrak A}}%
        \end{prooftree}}\\%
      {\begin{prooftree}%
          \Hypo{\hphantom{\mathfrak A\land{}}\mathfrak C\rightarrow\overline{\mathfrak A}}\Infer[left label=$n)$]1{\mathfrak A\Land\mathfrak C\rightarrow\hphantom{\overline{\mathfrak A}}}%
        \end{prooftree}}&&&%
      {\begin{prooftree}%
          \Hypo{\mathfrak C\rightarrow(\mathfrak x)\,\mathfrak A(\mathfrak x)}\Infer[left label=$o)$]1{\mathfrak C\rightarrow\rlap{$\mathfrak A(n)$}\hphantom{(\mathfrak x)\,\mathfrak A(\mathfrak x)}}%
        \end{prooftree}}%
    \end{aligned}%
  \end{gathered}%
\]
Zu den Grundrelationen kann $\mathfrak A\rightarrow\mathfrak A$ hinzugenommen werden.

Dieses Ergebnis aus §~2 können wir jetzt ergänzen:%
\begin{enumerate}[label=\arabic*)]%
\item es kann auch die Schlußregel \ $\begin{prooftree}[center]\Hypo{A(\mathfrak a)}\Infer[left label=$p)$]1{A(n)}\end{prooftree}$ hinzugenommen
werden.%
\end{enumerate}

Der Beweis wird wieder durch eine \emph{transfinite Prämissen\bs
induktion} geführt. Ist $A(\mathfrak a)$ herleitbar in~$\rmN$ und ist die
letzte Schlußregel dieser Herleitung nicht%
\[%
  \begin{prooftree}%
    \Hypo{A(1)}\Hypo{\cdots}\Hypo{A(n)}\Hypo{\cdots}\Infer4{A(\mathfrak a)}%
  \end{prooftree}%
\]%
so hat die Prämisse die Form~$A'(\mathfrak a)$. Nehmen wir als Induk\bs
tionsvoraussetzung an, daß für jede Prämisse $A'(\mathfrak a)$ auch $A'(n)$
herleitbar ist, so folgt sofort $A(n)$.\bruch{12}

\begin{enumerate}[resume*]%
\item Zu den Grundrelationen darf $\overline{\overline{\mathfrak A}}\rightarrow\mathfrak A$ hinzugenommen
werden.%
\end{enumerate}

Für jede Primformel~$\mathfrak P$ gilt nämlich stets ${}\rightarrow\mathfrak P$ \emph{oder} $\mathfrak P\rightarrow{}$.
Wegen \smash{$\begin{prooftree}[center]\Hypo{\hphantom{\overline{\overline{\mathfrak P}}}\rightarrow\mathfrak P}\Infer1{\overline{\overline{\mathfrak P}}\rightarrow\mathfrak P}\end{prooftree}$} \quad $ \begin{prooftree}[center] \Hypo{\mathfrak P\rightarrow\hphantom{\overline{\mathfrak P}}}\Infer1{\hphantom{\overline{\overline{\mathfrak P}}}\rightarrow\overline{\mathfrak P}}\Infer1{\overline{\overline{\mathfrak P}}\rightarrow\mathfrak P} \end{prooftree} $ \ 
ist also für jede Primformel stets $\overline{\overline{\mathfrak P}}\rightarrow\mathfrak P$  herleitbar. Hieraus
folgt allgemein die Herleitbarkeit von $\overline{\overline{\mathfrak A}}\rightarrow\mathfrak A$ (vergl.\ etwa
Hilbert-Bernays, \emph{Grundlagen der Mathematik~II}).

\begin{enumerate}[resume*]%
\item Es kann auch die \emph{vollständige Induktion}%
  \begin{prooftree*}%
    \Hypo{\mathfrak A(\mathfrak a)\rightarrow\mathfrak A(\mathfrak a')}\Infer[left label=$q)$]1{\mathfrak A(1)\rightarrow\mathfrak A(\mathfrak b)}%
  \end{prooftree*}%
zu den Schlußregeln hinzugenommen werden ohne
die Menge der herleitbaren Relationen zu vergrößern.%
\end{enumerate}

Ist nämlich $\mathfrak A(\mathfrak a)\rightarrow\mathfrak A(\mathfrak a')$ herleitbar, so auch die Relation
$\mathfrak A(n)\rightarrow\mathfrak A(n')$ für jede Zahl~$n$.

Für jede Zahl~$m$ folgt daraus durch $m$-malige Anwendung der
Schlußregel $k)$ sofort $\mathfrak A(1)\rightarrow\mathfrak A(m)$.

Wegen ${\begin{prooftree}[center] \Hypo{\mathfrak A(1)\rightarrow\mathfrak A(1)}\Hypo{\cdots}\Hypo{\mathfrak A(1)\rightarrow\mathfrak A(m)}\Hypo{\cdots}\Infer4{\mathfrak A(1)\rightarrow\mathfrak A(\mathfrak b)} \end{prooftree}}$
ist also auch $\mathfrak A(1)\rightarrow\mathfrak A(\mathfrak b)$ herleitbar.

Damit ist die Wf.\ der reinen Zahlentheorie bewiesen, da die
insgesamt zulässigen Schlußregeln einen Kalkül definieren,
der den klassischen Prädikatenkalkül ersichtlich enthält.\vfill

}

\ParallelRText{\renewcommand\red[1]{{\fontseries{b}\selectfont #1}}
\label{debut}\rightline{[P. LORENZEN]\qquad\qquad}\medskip

\begin{center}
  \uline{A proof of freedom from contradiction within the theory of partial order.}
\end{center}
  
% \looseness=-1
The dissertation of G. Gentzen contains a proof of freedom from contradiction
of elementary number theory without complete induction that relies
on the following basic thought: every derivable
sequent must also be derivable without detour, so that
during the derivation only those connectives are being introduced
that are absolutely\Deutsch{unbedingt} necessary, i.e.\ those\Deutsch{diejenigen} that are contained in the
sequent itself. In the proof of freedom from contradiction of number theory with complete induction, this basic thought
steps back with regard to others. I wish however to show in the following
that it alone suffices to obtain also this freedom from contradiction.

Without\label{beginning-2nd} knowledge of the dissertation of Gentzen, I have arrived at
this possibility on the basis of a semilattice-theoretic
question. This question is: how may a partially ordered set be embedded into an orthocomplemented complete
semilattice? In general, several such embeddings
are possible --~but among the possible embeddings one
is distinguished, i.e.\ the one which may be mapped
homomorphically into every other. The existence of this distinguished embedding will be proved in §~2.

In order to obtain from this in §~3 the sought-after proof of freedom from contradiction,
now just\Deutsch{nur noch} a translation of the semilattice-theoretic
proof into the logistic language is necessary. For the calculus
that we consider, and to which the usual calculuses may be reduced, is contained in the distinguished embedding of the
partially ordered set of the number-theoretic prime formulas.\bruch{2}\medskip

\uline{§\enspace1.} A set~$\rmM$ is called \red{partially ordered} if a binary relation~$\leqslant$ is defined in~$\rmM$ so that for the elements
$a,b,\dots$ of~$\rmM$ holds:
\[\begin{gathered}
  a\leqslant a\\
  a\leqslant b,\;b\leqslant c\quad\Rightarrow\quad a\leqslant c\text.
\end{gathered}\]
If $a\leqslant b$ and $b\leqslant a$ holds, then we write $a\equiv b$.

If $a\leqslant x$ holds for every~$x\in\rmM$, then we write $a\leqslant{}$. We write
as well ${}\leqslant a$ if $x\leqslant a$ holds for every~$x$. (${}\leqslant{}$ means
thus that $x\leqslant y$ holds for every~$x,y\in\rmM$.)

A partially ordered set~$\rmM$ is called \red{semilattice} if to every
$a,b\in\rmM$ there is a~$c\in\rmM$ so that for every~$x\in\rmM$ holds
\[x\leqslant a,\;x\leqslant b\quad\Longleftrightarrow\quad x\leqslant c\text.\]
$c$ is called the \blue{conjunction} of~$a$ and~$b$: \ $c\equiv a\land b$.

A semilattice~$\rmM$ is called \red{orthocomplemented}\label{def:orthocomplemented} if to every~$a\in\rmM$ there is
a~$b\in\rmM$ so that for every~$x\in\rmM$ holds
\[a\land x\leqslant{}\quad\Longleftrightarrow\quad x\leqslant b\text.\]
$b$~is called the \blue{orthocomplement} of~$a$: \ $b\equiv\neg{a}$.

A semilattice~$\rmM$ is called \red{$\omega$-complete} if to every countable sequence $\itM=a_1,a_2,\dots$ in~$\rmM$ there is a~$c\in\rmM$ so that for every~$x\in\rmM$ holds:
\[(\text{for every~$n$: }x\leqslant a_n)\quad\Longleftrightarrow\quad x\leqslant c\text.\]
% \looseness=-1 $c$~is called the \blue{conjunction of the elements of}~$\itM$: \ $c\equiv\bigland_na_n\equiv\bigland_{\itM}{}$.

% \looseness=1
If $\rmM$~and~$\rmM'$ are partially ordered sets, then $\rmM$~is called \red{a part} of~$\rmM'$
if~$\rmM$ is a subset of~$\rmM'$ and for every~$a,b\in\rmM$ \ $a\leqslant b$ holds
in~$\rmM'$ exactly if $a\leqslant b$ holds in~$\rmM$.

If $\rmM$~and~$\rmM'$ are partially ordered sets, we understand by
a \red{mapping} of~$\rmM$ into~$\rmM'$ an assignment that to every~$a\in\rmM$ assigns
an~$a'\in\rmM'$ so that
\[a\equiv b\quad\Rightarrow\quad a'\equiv b'\text.\bruch{3}\]

If $\rmM$~and~$\rmM'$ are orthocomplemented $\omega$-complete semilattices, we understand by a \red{homomorphism} of~$\rmM$
into~$\rmM'$ a mapping~$\to$ of~$\rmM$ into~$\rmM'$, so that for every~$a,b\in\rmM$
and~$a',b'\in\rmM'$ with $a\to a'$ and $b\to b'$ holds:
\[\begin{aligned}
  a\land b&\to a'\land b'\\
  \neg a&\to\neg{a'}\text.
\end{aligned}\]
Moreover, for every sequence~$\itM=a_1,a_2,\dots$ in~$\rmM$ and
$\itM'=a'_1,a'_2,\dots$ in~$\rmM'$ with $a_n\to a'_n$ is to hold:
\[\bigland_{\itM}{}\to\bigland_{\itM'}{}\text.\]

We want to prove now that to every partially ordered
set~$\rmP$ there is an orthocomplemented $\omega$-complete semilattice~$\rmK$ so that
\begin{itemize}
\item[1)]$\rmP$~is a part of~$\rmK$, %[i.e.\ $\rmK$~contains~$\rmP$ as part,]
\item[2)]$\rmK$~may be mapped homomorphically into every orthocomplemented $\omega$-com\-plete semilattice that contains $\rmP$~as part.
\end{itemize}
% \looseness=1
If $\rmK'$~were a further orthocomplemented $\omega$-complete
semilattice that fulfils conditions 1)~and~2), then
there would be an assignment by which $\rmK$~would be mapped homomorphically into~$\rmK'$ and $\rmK'$~into~$\rmK$, i.e.\ $\rmK$~and~$\rmK'$ would be \red{isomorphic}. $\rmK$~is thus
determined uniquely up to isomorphism by conditions 1)~and~2). We call $\rmK$~\red{the distinguished orthocomplemented $\omega$-complete semi\-lattice over~$\rmP$}.\medskip

\uline{§\enspace2.} \uline{Theorem}: \blue{There is over every partially ordered set the distinguished orthocomplemented\/ $\omega$-complete semilattice.}

We construct for the partially ordered set~$\rmP$ a set~$\rmK$
in the following way:
\begin{enumerate}[label=\arabic*)]
\item %\looseness=1
  Let $\rmK$~contain the elements of~$\rmP$. (These we call the
\red{prime elements} of~$\rmK$.)\bruch{4}
\item %\looseness=1
  Let $\rmK$~contain with finitely many elements~$a_1,a_2,\dots,a_n$
also the \red{combination} formed out of these as element.
(These we designate by~$a_1\land a_2\land\cdots\land a_n$.)
\item Let $\rmK$~contain with every element~$a$ also an element~$\neg a$.
\item Let $\rmK$~contain with every countable sequence~$\itM$ also an
element~$\bigland_{\itM}{}$.
\end{enumerate}

Every element of~$\rmK$ may thus be written uniquely as combination~$a_1\land a_2\land\cdots\land a_n$ of prime elements and elements of the form~$\neg a$
or~$\bigland_{\itM}{}$.

We define a relation~$\leqslant$ in~$\rmK$ in the following way:
\begin{itemize}
\item[1)]For prime elements~$p,q$ let $p\leqslant q$ hold in~$\rmK$ if~$p\leqslant q$ holds
  in~$\rmP$. (These relations we call the basic relations.)
\item[2)]Every relation~$\leqslant$ that may be derived from the
basic relations by the aid of the following rules is to hold in~$\rmK$:
\[
  \begin{gathered}
    \begin{aligned}[t]
      &{\begin{prooftree}
          \Hypo{c\leqslant a}\Hypo{c\leqslant b}\Infer[left label=$a)$]2{c\leqslant a\land b}
        \end{prooftree}}\\
      &{\begin{prooftree}
          \Hypo{a\land c\leqslant\hphantom a\vphantom b}\Infer[left label=$b)$]1{\vphantom{b}\hphantom{a\land{}}c\leqslant\neg a}
        \end{prooftree}}\\
      &{\label{rule-c}\begin{prooftree}
          \Hypo{\vphantom{\neg{b}}c\leqslant a_1}\Hypo{\cdots}\Hypo{c\leqslant a_n}\Hypo{\cdots}\Infer[left label=$c)$]4{c\leqslant\bigland_{\itM}}
        \end{prooftree}}\\
    \end{aligned}\quad
    \begin{aligned}[t]
      &{\begin{prooftree}
          \Hypo{\hphantom{b\land{}}a\leqslant c}\Infer[left label=$d{})$]1{a\land b\leqslant c}
        \end{prooftree}}\\
      &{\begin{prooftree}
          \Hypo{\hphantom{\neg b\land{}}a\leqslant b}\Infer[left label=$e)$]1{a\land\neg b\leqslant c}
        \end{prooftree}}\\
      &{\begin{prooftree}
          \Hypo{a_n\land b\leqslant c}\Infer[left label=$f)$]1{{\bigland_{\itM}}\land b\leqslant c}
        \end{prooftree}}
    \end{aligned}\\
    (\itM=a_1,a_2,\dots)\\
    {\begin{prooftree}
        \Hypo{a\land a\land b\leqslant c}\Infer[left label=$g)$]1[\bruch{5}]{\hphantom{a\land{}}a\land b\leqslant c}
      \end{prooftree}}
  \end{gathered}
\]
\end{itemize}
\looseness=1
We call the relations above the line the \red{premisses} of
the relation below the line.

We have now to show first that~$\rmK$
is an orthocomplemented $\omega$-complete semilattice w.r.t.\ the relation~$\leqslant$. For this we must prove
\[
  \begin{aligned}
    \alpha)&&a\leqslant a&\\
    \beta)&&a\leqslant b,\;b\leqslant c\quad&\Rightarrow\quad a\leqslant c\label{rule-beta}\\
    \gamma)&&c\leqslant a\land b\quad&\Rightarrow\quad c\leqslant a\\
    \delta)&&c\leqslant\neg a\quad&\Rightarrow\quad a\land c\leqslant{}\\
    \varepsilon)&&c\leqslant\bigland_{\itM}\quad&\Rightarrow\quad c\leqslant a_n&&(\itM=a_1,a_2,\dots)
  \end{aligned}
\]
These properties together with~$a)$, $b)$, and~$c)$ express in fact that~$\rmK$ is an orthocomplemented $\omega$-complete semilattice.

$\alpha)$ holds for prime elements. If $\alpha)$ holds for~$a$ and~$b$, then also for~$a\land b$ because of
\[
  \begin{prooftree}
  \Hypo{\hphantom{b\land{}}a\leqslant a}\Infer1{a\land b\leqslant a}\Hypo{\hphantom{a\land{}}b\leqslant b}\Infer1{a\land b\leqslant b}\Infer2{a\land b\leqslant a\land b}
\end{prooftree}
\]

If $\alpha)$ holds for every~$a_n\in\itM$, then also for~$\bigland_{\itM}$ because of
\[
  \begin{prooftree}
    \Hypo{a_1\leqslant a_1}\Infer1{{\bigland_{\itM}}\leqslant a_1}\Hypo{\cdots}\Infer[no rule]1{\vrule height5pt width0pt}%\rewrite{\raise2ex\box\treebox}
    \Hypo{a_n\leqslant a_n}\Infer1{{\bigland_{\itM}}\leqslant a_n}\Hypo{\cdots}\Infer[no rule]1{\vrule height5pt width0pt}%\rewrite{\raise2ex\box\treebox}
    \Infer4{{\bigland_{\itM}}\leqslant\bigland_{\itM}}
  \end{prooftree}
\]

If $\alpha)$ holds for~$a$, then also for~$\neg a$, because of
\[
  \begin{prooftree}
    \Hypo{\hphantom{\neg a\land{}}a\leqslant a}\Infer1{a\land\neg a\leqslant\hphantom{a}}\Infer1{\hphantom{a\land{}}\neg a\leqslant\neg a}
  \end{prooftree}
\]
Hereby $\alpha)$ is proved in general.\bruch{6}

As $\beta)$ is the most difficult to prove, we take first $\gamma)$.

In order to prove $\gamma)$, we have to show that if~$c\leqslant a\land b$
is derivable, then also $c\leqslant a$ must always be derivable.

We lead the proof indirectly\label{reductio} by a \blue{transfinite induction}. Let $c\leqslant a\land b$~be derivable, but not~$c\leqslant a$. Then the last
step of the derivation of $c\leqslant a\land b$ cannot be
${\begin{prooftree}[center]\Hypo{c\leqslant a}\Hypo{c\leqslant b}\Infer2{c\leqslant a\land b}\end{prooftree}}$, likewise\Deutsch{ebenfalls} not ${\begin{prooftree}[center]\Hypo{\hphantom{\neg{c_2}\land{}}c_1\leqslant\rlap{$c_2$}\hphantom{a\land b}}\Infer1{c_1\land\neg{c_2}\leqslant a\land b}\end{prooftree}}$ ($c=c_1\land\neg{c_2}$),
as then ${\begin{prooftree}[center]\Hypo{\hphantom{\neg{c_2}\land{}}c_1\leqslant c_2}\Infer1{c_1\land\neg{c_2}\leqslant\rlap{$a$}\hphantom{c_2}}\end{prooftree}}$ would be derivable at once.

For the last step remain only the possibilities
\[
  \begin{gathered}
    {\begin{prooftree}
        \Hypo{\hphantom{c_2\land{}}c_1\leqslant a\land b}\Infer1{c_1\land c_2\leqslant a\land b}
      \end{prooftree}}\qquad\qquad
      {\begin{prooftree}
        \Hypo{c_1\land c_1\land c_2\leqslant a\land b}\Infer1[$(c=c_1\land c_2)$]{\hphantom{c_1\land{}}c_1\land c_2\leqslant a\land b}
      \end{prooftree}}\\
      {\begin{prooftree}
          \Hypo{c_1\land c'\leqslant a\land b}\Hypo{\cdots}\Hypo{c_n\land c'\leqslant a\land b}\Hypo{\cdots}\Infer4[$\left(
        \begin{aligned}
          \itM&=a_1,a_2,\dots\\
          c&=\smash{\bigland_{\itM}}\land c'
        \end{aligned}
\right)$]{{\bigland_{\itM}}\land c'\leqslant a\land b}
        \end{prooftree}}
  \end{gathered}
\]

% \looseness=1
Here must now~$c_1\leqslant a$ resp.\ $c_1\land c_1\land c_2\leqslant a$ resp.\ for at least one~$n$ \ $c_n\land c'\leqslant a$ not be derivable, as otherwise at once
$c\leqslant a$ would be derivable. In the derivation of~$c\leqslant a\land b$ the claim~$\gamma)$ would thus
already be false for a premiss. If in
the derivation of a relation I go over to a premiss, of
this again to a premiss, etc., then I am after finitely
many steps at a basic relation. We would thus obtain
a basic relation for which the claim~$\gamma)$ would be false.
But as this is impossible, $\gamma)$~is thereby proved.

% \looseness=1
We call the induction that we have undertaken here
a \red{premiss induction}.

By the aid of premiss inductions, the proof for~$\delta)$ and~$\varepsilon)$ proceeds just as simply as for~$\gamma)$, so that I am not going into this any further\Deutsch{hierauf eingehe}.\bruch{7}

It remains only to show in addition $\beta)$. Instead of this we prove the
stronger claim
\[\zeta)\quad a\leqslant b,\; b\land b\land\cdots\land b\land c\leqslant d\quad \Rightarrow\quad a\land c\leqslant d\]
in order to be able to apply premiss inductions hereupon.

% \looseness=1
Let first~$b$, $c$~and~$d$ be prime elements. Then $\zeta)$~holds
for every basic relation~$a\leqslant b$. We assume as \blue{induction hypothesis} that $\zeta)$~holds \blue{for every premiss of}~$a\leqslant b$.

% \looseness=1
As $b$~is a prime element, the last step of the derivation of~$a\leqslant b$ can only be:
\[
  \begin{gathered}
    {\begin{prooftree}
        \Hypo{\hphantom{a_2\land{}}a_1\leqslant b}\Infer%[left label=(1)]
        1{a_1\land a_2\leqslant b}
      \end{prooftree}}\qquad
    {\begin{prooftree}
        \Hypo{a_1\land a_1\land a_2\leqslant b}\Infer%[left label=(2)]
        1[$(a=a_1\land a_2)$]{\hphantom{a_1\land{}}a_1\land a_2\leqslant b}
      \end{prooftree}}\\
    {\begin{prooftree}
        \Hypo{\hphantom{a_2\land{}}a_1\leqslant a_2\vphantom{a'}}\Infer%[left label=(3)]
        1[$(a=a_1\land\neg{a_2})$]{a_1\land\neg{a_2}\leqslant \rlap{$b$}\hphantom{a_2}\vphantom{\bigland_{\itM}}}
      \end{prooftree}}\qquad
    {\begin{prooftree}
        \Hypo{a_n\land a'\leqslant b}\Infer%[left label=(4)]
        1[$\left(
      \begin{aligned}
        \itM&=a_1,a_2,\dots\\
        a&=\smash{\bigland_{\itM}}\land a'
      \end{aligned}
    \right)$]{{\bigland_{\itM}}\land a'\leqslant b}
      \end{prooftree}}
  \end{gathered}
\]

According to the induction hypothesis, then $a_1\land c\leqslant d$ %~~(1)
resp.\ $a_1\land a_1\land a_2\land c\leqslant d$ %~~(2)
resp.\ $a_n\land a'\land c\leqslant d$ %~~(4)
is derivable.
In every case $a\land c\leqslant d$ is at once derivable, as well from $a_1\leqslant a_2$ 
because of
\[\begin{prooftree}
  \Hypo{\hphantom{a_2\land{}}a_1\leqslant a_2}\Infer1{a_1\land\neg{a_2}\leqslant\rlap{$d$}\hphantom{a_2}}\Infer1{\hphantom{a_1\land a_2}\llap{$a\land c$}\leqslant\rlap{$d$}\hphantom{a_2}}
\end{prooftree}
\]

Thereby $\zeta)$~is proved for prime elements~$b$, $c$~and~$d$.

Now let only $b$~still be a prime element. Then $\zeta)$~holds thus for arbitrary~$a$ and prime elements~$c$,~$d$. A \blue{premiss induction} results now in $\zeta)$~holding for every relation $b\land b\land\cdots\land b\land c\leqslant d$. Every premiss of $b\land b\land\cdots\land b\land c\leqslant d$ has in fact again
the form $b\land\cdots\land b\land c\leqslant d$. Thereby $\zeta)$ is proved in general for \blue{prime elements}~$b$.

% \looseness=1
If $\zeta)$ holds for elements~$b_1$ and~$b_2$, then obviously also for $b_1\land b_2$.
If $\zeta)$~holds for every~$b_n\in\itM$, then also for~$b=\bigland_{\itM}$. (Proof by
premiss induction: ${\bigland_{\itM}}\land{\bigland_{\itM}}\land\dots\land{\bigland_{\itM}}\land c\leqslant d$ can have the following
premiss: $b_n\land{\bigland_{\itM}}\land\cdots\land{\bigland_{\itM}}\land c\leqslant d$. According to induction hypo\bruch{8}%
thesis holds then $a\leqslant{\bigland_{\itM}},\;b_n\land{\bigland_{\itM}}\land\cdots\land{\bigland_{\itM}}\land c\leqslant d\;\Rightarrow\;b_n\land a\land c\leqslant d$.
But as $\zeta)$ is also assumed for $b=b_n$, and because of
\[a\leqslant{\bigland_{\itM}}\quad\Rightarrow\quad a\leqslant b_n\text,\]
also $a\leqslant{\bigland_{\itM}},\;b_n\land a\land c\leqslant d\;\Rightarrow\;a\land a\land c\leqslant d$ holds.
But from $a\land a\land c\leqslant d$ may be derived $a\land c\leqslant d$. Every other
premiss of ${\bigland_{\itM}}\land{\bigland_{\itM}}\land\cdots\land{\bigland_{\itM}}\land c\leqslant d$ is trivial.)

If $\zeta)$~holds for~$b$, then also for~$\neg b$. (Proof by premiss induction: $\neg b\land\neg b\land\cdots\land\neg b\land c\leqslant d$ can have the following premiss:
$\neg b\land\cdots\land\neg b\land c\leqslant b$. Then holds according to induction hypothesis
\[a\leqslant\neg b,\;\neg b\land\cdots\land\neg b\land c\leqslant b\quad\Rightarrow\quad a\land c\leqslant b\text.\]
As $\zeta)$ is also assumed for~$b$, also holds
\[a\land c\leqslant b,\;a\land b\leqslant d\quad\Rightarrow\quad a\land a\land c\leqslant d\text.\]
Thus holds also $a\leqslant\neg b,\;\neg b\land\cdots\land\neg b\land c\leqslant b\;\Rightarrow\;a\land c\leqslant d$
because of $a\leqslant\neg b\;\Rightarrow\;a\land b\leqslant d$.
Every other premiss is again trivial.)

Thus $\zeta)$~is valid in general. This proves that $\rmK$~is an
orthocomplemented $\omega$-complete semilattice.

% \looseness=1
$\rmP$~is a part of~$\rmK$, as
\[
  p\leqslant q\text{ in }\rmP\quad\Longleftrightarrow\quad p\leqslant q\text{ in }\rmK
\]
holds. We have for this to convince ourselves that no relation
$p\leqslant q$ is derivable in~$\rmK$ that is not already holding in~$\rmP$. But this goes without saying\Deutsch{selbstverständlich sein}, as none of the rules except~$g)$ actually\Deutsch{überhaupt} yields relations~$p\leqslant q$ below the line. A derivation of a relation~$p\leqslant q$ can thus use
only the rules~$d)$ and~$g)$. But with these only the basic relations are derivable.

\looseness=1
For the proof of our theorem, it remains now in addition to show that
$\rmK$~may be mapped homomorphically into every other orthocomplemented $\omega$-complete
semilattice~$\rmK'$ that contains~$\rmP$ as part\bruch{9}. This mapping we define by
\begin{enumerate}[label=\arabic*)]
\item for prime elements~$p$ holds~$p\to p$,
\item moreover is to hold
  \[
    \begin{aligned}
      a\to a',\;b\to b'\quad&\Rightarrow\quad a\land b\to a'\land b'\\
      a\to a'\quad&\Rightarrow\quad\neg a\to\neg{a'}\\
      a_n\to a'_n\quad&\Rightarrow\quad\bigland_{\itM}\to\bigland_{\itM'}\quad
      \left(
        \begin{aligned}
          \itM&=a_1,a_2,\dots\\
          \itM'&=a'_1,a'_2,\dots
        \end{aligned}
      \right)
    \end{aligned}
  \]
\end{enumerate}
Hereby obviously a homomorphism is being defined, for
with~$a\to a'$ and~$b\to b'$ always holds~$a\leqslant b\Rightarrow a'\leqslant b'$.

% \looseness=1
Every derivation of $a\leqslant b$ proves in fact at once also $a'\leqslant b'$,
as the derivation steps~$a)$--$g)$ are always correct in every orthocomplemented $\omega$-complete semilattice.\medskip

\uline{§\enspace3.} In order to be able to prove the freedom from contradiction of elementary number theory with complete induction
from the theorem proved in~§~2, we use the following formalisation.
We take as prime formulas the signs for number-theoretic
predicates~$\rmA(\dots)$, $\rmB(\dots)$,~\ldots\ with the numbers $1,1',1'',\dots$
as arguments, e.g.\ $1 = 1''$, $1+1=1'$.

These prime formulas~$\mathfrak{P},\mathfrak{Q},\dots$ form a partially ordered set
if we set $\mathfrak{P}\rightarrow\mathfrak{Q}$ in case the predicate~$\mathfrak{P}$
implies the predicate~$\mathfrak{Q}$. To the basic relations $\mathfrak{P}\rightarrow\mathfrak{Q}$ we are also adding the relations of the form ${}\rightarrow\mathfrak{P}$, $\mathfrak{P}\rightarrow{}$, ${}\rightarrow{}$,
as far as they are correct in terms of content.

Over this partially ordered set~$\rmP$ of the prime formulas, we construct now as in §~2 the distinguished orthocomplemented $\omega$-complete semilattice. We use for this
the logistic signs, thus~$\rightarrow$ instead of~$\leqslant$, $\Land$ instead of~$\land$.

\looseness=1
To the formulas belong thus the prime formulas, with~$\mathfrak{A}$ and~$\mathfrak{B}$
also $\mathfrak{A}\Land\mathfrak{B}$, with~$\mathfrak{A}$ also~$\overline{\mathfrak{A}}$. We restrict the conjunction of countable\bruch{10}
sequences to the sequences of the form $\mathfrak{A}(1),\mathfrak{A}(1'),\dots$.
We designate this conjunction by $(\mathfrak{x})\,\mathfrak{A}(\mathfrak{x})$.

Moreover, we introduce in addition free variables\label{free} $\mathfrak{a}$ = $a$, $b$,~\ldots\ by
the following rule of inference:
\[
  \pbox{215pt}{if $A(1),A(1'),\dots$ are derivable relations, then\\
    $A(\mathfrak{a})$ is also to be derivable.}
\]

By this the proofs of §~2 are only modified unessentially\Deutsch{unwesentlich}. We obtain overall\Deutsch{insgesamt} a calculus~$\rmN$\label{N} with the following rules of inference
\[
  \begin{gathered}
    \begin{aligned}[t]
      &{\begin{prooftree}
          \Hypo{\mathfrak{C}\rightarrow\mathfrak{A}}\Hypo{\mathfrak{C}\rightarrow\mathfrak{B}}\Infer[left label=$a)$]2{\mathfrak{C}\rightarrow\mathfrak{A}\Land\mathfrak{B}}
        \end{prooftree}}\\
      &{\begin{prooftree}
          \Hypo{\mathfrak{A}\Land \mathfrak{C}\rightarrow\hphantom{\overline{\mathfrak{A}}}}\Infer[left label=$b)$]1{\hphantom{\mathfrak{A}\Land{}}\mathfrak{C}\rightarrow\overline{\mathfrak{A}}}
        \end{prooftree}}\\
      &{\begin{prooftree}
          \Hypo{\mathfrak{C}\rightarrow\mathfrak{A}(1)}\Hypo{\cdots}\Hypo{\mathfrak{C}\rightarrow\mathfrak{A}(n)}\Hypo{\cdots}\Infer[left label=$c)$]4{\mathfrak{C}\rightarrow(\mathfrak{x})\,\mathfrak{A}(\mathfrak{x})}
        \end{prooftree}}\\
    \end{aligned}\quad
    \begin{aligned}[t]
      &{\begin{prooftree}
          \Hypo{\hphantom{\mathfrak{B}\Land{}}\mathfrak{A}\rightarrow \mathfrak{C}}\Infer[left label=$d{})$]1{\mathfrak{A}\Land\mathfrak{B}\rightarrow \mathfrak{C}}
        \end{prooftree}}\\
      &{\begin{prooftree}
          \Hypo{\hphantom{\overline{\mathfrak{B}}\Land{}}\mathfrak{A}\rightarrow\mathfrak{B}}\Infer[left label=$e)$]1{\mathfrak{A}\Land\overline{\mathfrak{B}}\rightarrow\mathfrak{C}}
        \end{prooftree}}\\
      &{\begin{prooftree}
          \Hypo{\hphantom{(\mathfrak{x})\,\mathfrak{A}(\mathfrak{x})}\llap{$\mathfrak{A}(n)$}\Land\mathfrak{B}\rightarrow \mathfrak{C}}\Infer[left label=$f)$]1{(\mathfrak{x})\,\mathfrak{A}(\mathfrak{x})\Land\mathfrak{B}\rightarrow \mathfrak{C}}
        \end{prooftree}}
    \end{aligned}\\
    {\begin{prooftree}
        \Hypo{\mathfrak{A}\Land\mathfrak{A}\Land\mathfrak{B}\rightarrow \mathfrak{C}}\Infer[left label=$g)$]1{\hphantom{\mathfrak{A}\Land{}}\mathfrak{A}\Land\mathfrak{B}\rightarrow \mathfrak{C}}
      \end{prooftree}}\\
    {\begin{prooftree}
        \Hypo{\mathfrak{A}\Land\mathfrak{B}\rightarrow\mathfrak{C}}\Infer[left label=$h)$]1{\mathfrak{B}\Land\mathfrak{A}\rightarrow\mathfrak{C}}
      \end{prooftree}}\qquad
        {\begin{prooftree}
        \Hypo{\mathfrak{A}\Land(\mathfrak{B}\Land\mathfrak{C})\rightarrow\mathfrak{D}}\Infer[left label=$i)$]1{(\mathfrak{A}\Land\mathfrak{B})\Land\mathfrak{C}\rightarrow\mathfrak{D}}
      \end{prooftree}}\\
    {\begin{prooftree}
          \Hypo{A(1)}\Hypo{\cdots}\Hypo{A(n)}\Hypo{\cdots}\Infer[left label=$j)$]4[\bruch{11}]{A(\mathfrak{a})}
        \end{prooftree}}
  \end{gathered}\label{rule-c-zt}
\]

The rules of inference~$h)$ and~$i)$ were dispensable\Deutsch{überflüssig} in §~2, as we
 have introduced there~$a\land b\land c\dots$ at once as sign for the combination of
$a,b,c,\dots$.

The proof in §~2 yields now the following result:
the calculus~$\rmN$ is consistent, e.g.\ the empty relation~${}\rightarrow{}$ is not derivable, as only the
relations correct in terms of content hold in~$\rmP$ and $\rmP$~is a part of~$\rmN$.
To the calculus~$\rmN$ the following rules of inference can be added without increasing the set the derivable relations:
\[
  \begin{gathered}
    {\begin{prooftree}
        \Hypo{\mathfrak{A}\rightarrow\mathfrak{B}}\Hypo{\mathfrak{B}\rightarrow\mathfrak{C}}\Infer[left label=$k)$]2{\mathfrak{A}\rightarrow\mathfrak{C}}
      \end{prooftree}}\label{rule-k}\\
    \begin{aligned}
      {\begin{prooftree}
          \Hypo{\mathfrak{C}\rightarrow\mathfrak{A}\Land\mathfrak{B}}\Infer[left label=$l)$]1{\mathfrak{C}\rightarrow\mathfrak{A}\hphantom{{}\Land\mathfrak{B}}}
        \end{prooftree}}&&\qquad&
      {\begin{prooftree}
          \Hypo{\mathfrak{C}\rightarrow\mathfrak{A}\Land\mathfrak{B}}\Infer[left label=$m)$]1{\mathfrak{C}\rightarrow\mathfrak{B}\hphantom{{}\Land\mathfrak{A}}}
        \end{prooftree}}\\
      {\begin{prooftree}
          \Hypo{\hphantom{\mathfrak{A}\land{}}\mathfrak{C}\rightarrow\overline{\mathfrak{A}}}\Infer[left label=$n)$]1{\mathfrak{A}\Land\mathfrak{C}\rightarrow\hphantom{\overline{\mathfrak{A}}}}
        \end{prooftree}}&&&
      {\begin{prooftree}
          \Hypo{\mathfrak{C}\rightarrow(\mathfrak{x})\,\mathfrak{A}(\mathfrak{x})}\Infer[left label=$o)$]1{\mathfrak{C}\rightarrow\rlap{$\mathfrak{A}(n)$}\hphantom{(\mathfrak{x})\,\mathfrak{A}(\mathfrak{x})}}
        \end{prooftree}}
    \end{aligned}
  \end{gathered}
\]
To the basic relations can be added $\mathfrak{A}\rightarrow\mathfrak{A}$.

This result from §~2 we can now complete\Deutsch{ergänzen}:
\begin{enumerate}[label=\arabic*)]
\item The rule of inference \ $\begin{prooftree}[center]\Hypo{A(\mathfrak{a})}\Infer[left label=$p)$]1{A(n)}\end{prooftree}$ can also be added.
\end{enumerate}

The proof is again being led by a \blue{transfinite premiss induction}. If $A(\mathfrak{a})$ is derivable in~$\rmN$ and if the
last rule of inference of this derivation is not
\[
  \begin{prooftree}
    \Hypo{A(1)}\Hypo{\cdots}\Hypo{A(n)}\Hypo{\cdots}\Infer4{A(\mathfrak{a})}
  \end{prooftree}
\]
then the premiss has the form~$A'(\mathfrak{a})$. If we assume as induction hypothesis\Deutsch{nehmen ... an} that for every premiss $A'(\mathfrak{a})$ also $A'(n)$
is derivable, then $A(n)$ follows at once.\bruch{12}

\begin{enumerate}[resume*]
\item To the basic relations may be added $\overline{\overline{\mathfrak{A}}}\rightarrow\mathfrak{A}$.
\end{enumerate}

For every prime formula~$\mathfrak{P}$ holds in fact always $\rightarrow\mathfrak{P}$ \blue{or} $\mathfrak{P}\rightarrow$.
Because of \ $\begin{prooftree}[center]\Hypo{\hphantom{\overline{\overline{\mathfrak{P}}}}\rightarrow\mathfrak{P}}\Infer1{\overline{\overline{\mathfrak{P}}}\rightarrow\mathfrak{P}}\end{prooftree}$ \quad $ \begin{prooftree}[center]\Hypo{\mathfrak{P}\rightarrow\hphantom{\overline{\mathfrak{P}}}}\Infer1{\hphantom{\overline{\overline{\mathfrak{P}}}}\rightarrow\overline{\mathfrak{P}}}\Infer1{\overline{\overline{\mathfrak{P}}}\rightarrow\mathfrak{P}}\end{prooftree}$, \ $\overline{\overline{\mathfrak{P}}}\rightarrow\mathfrak{P}$
is thus always derivable for every prime formula. From this
follows in general the derivability of $\overline{\overline{\mathfrak{A}}}\rightarrow\mathfrak{A}$ (cf.\ e.g.\ Hilbert-Bernays, \emph{Grundlagen der Mathematik~II}).

\begin{enumerate}[resume*]
\item The \blue{complete induction}\label{complete-induction}
  \[
    {\begin{prooftree}
      \Hypo{\mathfrak{A}(\mathfrak{a})\rightarrow\mathfrak{A}(\mathfrak{a}')}\Infer[left label=$q)$]1{\mathfrak{A}(1)\rightarrow\mathfrak{A}(\mathfrak{b})}
    \end{prooftree}}
  \]
can also be added to the rules of inference without increasing
the set the derivable relations.
\end{enumerate}

In fact, if $\mathfrak{A}(\mathfrak{a})\rightarrow\mathfrak{A}(\mathfrak{a}')$ is derivable, then also the relation
$\mathfrak{A}(n)\rightarrow\mathfrak{A}(n')$ for every number~$n$.

For every number~$m$ follows therefrom at once $\mathfrak{A}(1)\rightarrow\mathfrak{A}(m)$ by $m$-fold application of the
rule of inference~$k)$.

Because of ${\begin{prooftree}[center]\Hypo{\mathfrak{A}(1)\rightarrow\mathfrak{A}(1)}\Hypo{\cdots}\Hypo{\mathfrak{A}(1)\rightarrow\mathfrak{A}(m)}\Hypo{\cdots}\Infer4{\mathfrak{A}(1)\rightarrow\mathfrak{A}(\mathfrak{b})}\end{prooftree}}$
also $\mathfrak{A}(1)\rightarrow\mathfrak{A}(\mathfrak{b})$ is thus derivable.

Thereby the freedom from contradiction of the elementary number theory is proved, as the
overall\Deutsch{insgesamt} admissible rules of inference define a calculus
that obviously contains the classical calculus of predicates.\label{fin}\vfill

}
\end{Parallel}

\end{document}